\def\versiondate{30 March 2000}
\input math.macros
\input Ref.macros

\checkdefinedreferencetrue
\continuousfigurenumberingtrue
\theoremcountingtrue
\sectionnumberstrue
\forwardreferencetrue
\tocgenerationtrue
\citationgenerationtrue
\nobracketcittrue
\hyperstrue
\initialeqmacro

\bibsty{myapalike}
\input\jobname.lbl 

\def\vertex{V}
\def\edge{E}
\def\Aut{{\rm Aut}}
\def\\{\backslash}
\def\bd{\partial}
\def\iso{\iota}
\def\bde{\bd_\edge}
\def\bdv{\bd_\vertex}
\def\bdiv{\bd_\vertex^{\rm int}}  
\def\isoe{\iso_\edge}
\def\isov{\iso_\vertex}

\def\bp{o}
\def\path{{\cal P}}
\def\mod#1{|#1|_*}   

\def\pcb{p_c^{\rm bond}}
\def\pcs{p_c^{\rm site}}
\def\pub{p_u^{\rm bond}}
\def\pus{p_u^{\rm site}}
\def\gr{{\rm gr}}

\def\edges{E}
\def\ev#1{{\cal #1}}
\def\gdeg{d}  
\def\tdeg{b}  
\def\gp{\Gamma}  
\def\gpe{\gamma}  
\def\gh{G}  
\def\br{{\rm br}}
\def\CAP{{\rm cap}}
\def\t #1 {T_{#1}}   
\def\tu #1 {T^{#1}}    
\def\rc{{\ss RC}}  
\def\frc{{\ss FRC}}  
\def\wrc{{\ss WRC}}  
\def\fpt{{\ss FPt}}   
\def\ppt{{\ss Pt}}    
\def\wpt{{\ss WPt}}    

\def\RWl{${\ss RW}_\lambda$\ }
\def\degn{\deg_-}
\def\en{{\cal E}}   
\def\ncomp(#1){\|#1\|}  
\def\ising{{\ss Ising}}
\def\rsign{J}   
\def\SG{{\ss SpGl}}  
\def\cpr{\xi}    
\def\thetaP{\theta(\P)}  
\def\ance{\isoe^*} 
\def\ancv{\isov^*} 

\def\BLPSgip{\ref b.BLPSgip/\def\BLPSgip{\htmllocref{BLPSgip}{[BLPS99]}}}
\def\ACCN{\ref b.ACCN/\def\ACCN{\htmllocref{ACCN}{[ACCN]}}}

\def\firstheader{\eightpoint\ss\underbar{\raise2pt\line
    {{\it J.\ Math.\ Phys.} {\bf 41} (2000), 1099--1126.
    \hfil Version of \versiondate}}}

\beginniceheadline

\vglue20pt

\title{Phase Transitions on Nonamenable Graphs}

\author{Russell Lyons}

\abstract{We survey known results about phase transitions in various models
of statistical physics when the underlying space is a nonamenable graph.
Most attention is devoted to transitive graphs and trees.
}

\bottom{Primary
82B26
. Secondary
60B99, 
60F60, 
60J15, 
60K35, 
82B20, 
82B43, 
82C22, 
82C43
.}
{Percolation, amenable, unimodular, transitive graph, Ising, Potts, random
cluster, trees, contact process, random walk.}
{Research partially supported by
NSF grant DMS-9802663.}

\articletoc

\bsection {Introduction}{s.intro}

We shall give a summary of some of the main results known about phase
transitions on nonamenable graphs. All terms will be defined as needed
beginning in \ref s.back/.
Among the graphs we consider, we pay special attention to transitive
graphs and trees (regular or not), as
these are the cases that arise most naturally. Both of these classes of
graphs also have some feature that permits a satisfying analysis to be
performed (or to be conjectured):
transitive graphs look the same from each vertex, while trees lack cycles.
Certain phenomena are known to occur for all transitive nonamenable graphs,
others are conjectured to hold for all transitive nonamenable graphs, while
still others depend on different aspects of the graph. The subject is quite
rich because of the interplay between probabilistic models and geometry.
In particular, there is a greater variety of probabilistic behavior
possible on nonamenable graphs than on amenable graphs.
The area is developing vigorously, but a great deal remains to be discovered.
A number of parallels among different processes will be evident to the
reader, and consequently, a number of questions will suggest themselves.
We have, however, omitted all discussion of critical exponents.

The models we consider all involve a parameter.
Changing the parameter leads to qualitative changes of behavior. When
such a change occurs, we shall say that there is a {\bf phase transition}.
(Note:
in some publications, a phase transition is said to occur for a {\it
fixed\/}
parameter value when there is more than one Gibbs measure at that value. By
contrast, our term is not precisely defined.)
There is usually at least one {\bf critical value} for the parameter, i.e.,
a value separating two intervals of the parameter where there are
different qualitative behaviors on each side of the critical value. For the
most basic phase transitions, those that usually occur on amenable graphs,
\ref b.Hag:mrf/ showed that a phase transition occurs simultaneously in all
or none of the following models on any given graph, assuming only that the
graph has bounded degree: bond percolation, site percolation, the Ising
model, the Widom-Rowlinson model, and the beach model. However, what makes
nonamenable graphs truly distinctive is often the presence of a second
critical value that does not occur on amenable graphs. The extent to which
such behavior is understood varies widely from model to model and from
graph to graph.

There are various probabilistic characterizations known of
nonamenability. The first
such result was proved for the most basic probabilistic process, namely,
random walk, in the thesis of
\refbmulti{Kesten:amenA,Kesten:amenB}.
He showed
that a countable group $\gp$ is amenable iff the spectral radius is $1$ for
some (or every)
symmetric group-invariant random walk whose support generates
$\gp$.
The extension of Kesten's theorem to the setting of invariant random walks
on transitive graphs involves
unimodularity and has been studied by \ref b.SoardiWoess:rws/, \ref
b.Salvatori/, and \ref b.SCosteWoess/.
We shall return to random walks, now with a parameter, in \ref s.rwl/.

Due to lack of time, we were unable to survey results concerning branching
random walk, which has many similarities to results here and, indeed, has
inspired many of them.
We mention just one example: A group $\gp$ is amenable iff for some (or every)
symmetric group-invariant random walk with support generating $\gp$ and for
some (or every) tree $T$ with branching number larger than 1, the
associated $T$-indexed random walk on $\gp$ is recurrent. In particular,
this is the case for branching random walk corresponding to any
Galton-Watson branching process with mean larger than 1. See \ref
b.BP:tirw/ for definitions and a proof (which depends on Kesten's theorem
above).
This result inspired \ref g.pcpu/ by means of an intuitive analogy between
the range of a branching random walk and an infinite percolation cluster;
see the proof of Thm.~4 in \ref b.pyond/ for a direct relationship between
branching random walk and percolation.

We now give a somewhat more detailed
preview of some of the results to be surveyed. For ordinary
Bernoulli percolation on transitive amenable graphs, it is well known that when
there is a.s.\ an infinite cluster, then there is a.s.\ a unique infinite
cluster. This is now known to fail in many cases of transitive nonamenable
graphs, and has been conjectured to fail in all transitive nonamenable
graphs. Moreover, it is known that the uniqueness and nonuniqueness phases,
if not empty, determine single intervals of the parameter. This leads to
the study of two critical parameters, the usual one at the top of the
regime of nonexistence of infinite clusters and a possibly new one at the
bottom of the uniqueness phase. It also leads to the study of the behavior
of the infinite clusters when there are infinitely many and how they merge
as the parameter is increased.

One of the important new tools for studying percolation is the
Mass-Transport Principle and its use in invariant percolation. This
provides some general results that allow one to manipulate the clusters of
Bernoulli percolation in a rather flexible fashion. In particular,
nonamenability turns out to be more of an asset than a liability, as it
provides for new thresholds that are trivial in the amenable case.

The Ising model is one of a natural family of models that includes
Bernoulli percolation.  Additional complications, such as boundary
conditions and the optional parameter of an external field, lead to
questions that do not arise for Bernoulli percolation. Sometimes, they can
be used to characterize exactly amenability. But
the number of different phase transitions that are possible for
Potts models and the related random cluster models is sufficiently great
that it has so far precluded the kind of unified picture that is at least
conjectured for percolation.

The contact process is now reasonably well understood on the euclidean
lattices $\Z^d$. However, some fundamental results there are still not
known in the more general setting of amenable transitive graphs. For
example, analogous to the number of infinite clusters in Bernoulli
percolation are the phases in the contact process of extinction, weak
survival, and strong survival. It might be that weak survival is impossible
iff the transitive graph is amenable. Some results in this direction are
known.

When we consider trees, we are often able to calculate precisely many
critical values for various processes, even for completely general trees
without any regularity. In almost all instances, these
critical values turn out to be functions of a single number associated to
the tree, its average branching number.


\bsection{Background on Graphs}{s.back}

The basic definitions of the terms pertaining to graphs are as follows.
Let $\gh = (\vertex, \edges)$ be an unoriented
graph with vertex set $\vertex$ and
symmetric edge set $\edge \subseteq \vertex \times \vertex$.
We write edges as $[x, y]$. If $x$ and $y$ are the endpoints of an edge, we
call them {\bf adjacent} or {\bf neighbors} and write $x \sim y$.
All graphs are assumed without
further comment to be connected, denumerable, and locally finite.
The only exception is that random subgraphs of a given graph may well be
disconnected.
Given $K \subset \vertex$, set $\bdv K := \{ y \notin K \st \exists x
\in K, \ x \sim y \}$ and
 $\bde K := \{[x,y] \in \edge \st x \in K, \,  y \notin K \}$.
Define the {\bf vertex-isoperimetric constant} of $\gh$ by
$$
\isov(\gh) := \inf \left\{ {|\bdv K| \over |K|} \st
  K \subset \vertex \hbox{ is finite and nonempty} \right\}\,,
$$
and  let the {\bf edge-isoperimetric constant} of $\gh$ be
$$
\isoe(\gh) := \inf \left\{ {|\bde K| \over |K|} \st
  K \subset \vertex \hbox{ is finite and nonempty} \right\}\,.
$$

A graph $\gh$ is
called {\bf amenable} if $\isoe(\gh) = 0$. If $\gh$ has
bounded degree, then this is equivalent to $\isov(\gh) = 0$. An {\bf
automorphism}
of $\gh$ is a bijection of $\vertex$ that induces a bijection of $\edges$.
The set of automorphisms of $\gh$ forms a group denoted $\Aut(\gh)$.
We say that a group $\gp
\subseteq \Aut(\gh)$ is {\bf transitive} or {\bf acts transitively} if
$\vertex$ has only one orbit under $\gp$, i.e., if for all $x, y \in
\vertex$, there is some $\gamma \in \gp$ such that $\gamma x = y$.
We say that $\gp$ is {\bf quasi-transitive} 
if $\gp$ splits $\vertex$ into finitely many orbits.
We call the graph $\gh$ itself {\bf (quasi-)transitive} if $\Aut(\gh)$ is.
Most results concerning
quasi-transitive graphs can be deduced from corresponding results for
transitive graphs or can be deduced in a similar fashion but with some
additional attention to details. For simplicity, we shall therefore ignore
quasi-transitive graphs in the sequel.
(The extension of results to quasi-transitive graphs is important, however.
Not only do they arise naturally, but they are crucial to the study of
planar transitive graphs.)

Let $\gp$ be a finitely generated group and $S$ a finite symmetric
generating set for
$\gp$.  The (right) {\bf Cayley graph} $\gh=\gh(\gp,S)$ of $\gp$
is the graph with vertex set $\vertex:=\gp$ and edge set
$\edges:=\bigl\{[v,vs]\st v\in\gp,\, s\in S\bigr\}$.
Note that $\gp$ acts transitively on
$\gh$ by the {\bf translations} $\gamma : x \mapsto \gamma x$.

A {\bf tree} is a graph without cycles or loops. A branching
process with one initial progenitor gives rise naturally to a random tree,
its genealogical tree. When the branching process is a Galton-Watson
process, we call the resulting tree a {\bf Galton-Watson tree}.

We now review the modular function.
Each compact group has a unique left-invariant Radon probability measure,
called Haar measure.  It is also
the unique right-invariant Radon probability measure.  A locally compact
group $\gp$ has a left-invariant $\sigma$-finite Radon measure $|\cbuldot|$;
it is unique up to a multiplicative constant.  For every $\gpe \in \gp$, the
measure $A \mapsto |A\gpe|$ is left invariant, whence there is a positive
number $m(\gpe)$ such that $|A\gpe| = m(\gpe) |A|$ for all measurable $A$.  The map
$\gpe \mapsto m(\gpe)$ is a homomorphism from $\gp$ to the multiplicative group
of the positive reals and is called the {\bf modular function} of $\gp$.
If $m(\gpe)=1$ for every $\gpe\in \gp$, then $\gp$ is called {\bf unimodular}.
In particular, this is the case if $\gp$ is countable, where Haar measure is
counting measure. See, e.g., \ref b.Royden/ for more on Haar measure.

We give the automorphism group $\Aut(\gh)$ of a
graph $\gh$ the topology of pointwise convergence.
By Corollary 6.2 of \ref b.BLPSgip/, if there is a transitive unimodular
closed subgroup of $\Aut(\gh)$, then $\Aut(\gh)$ is also unimodular. In
particular, this is the case if $\gh$ is the Cayley graph of a group $\gp$.
For this reason and for simplicity, we shall not generally consider
subgroups of $\Aut(\gh)$.
However, the reader may wish instead to concentrate on translation-invariant
measures on Cayley graphs, i.e., on the subgroup $\gp$ of automorphisms of
a Cayley graph $\gh$ of $\gp$.
We call a graph $\gh$ {\bf unimodular} if $\Aut(\gh)$ is.

The stabilizer
$$
S(x) := \{\gpe \in \Aut(\gh) \st \gpe x = x\}
$$
of any vertex $x$ is compact and so has finite Haar measure.
Note that if $\gpe u = y$, then $S(y) = \gpe S(u) \gpe^{-1}$, whence
$$
|S(y)| = |S(u) \gpe^{-1}| = m(\gpe)^{-1} |S(u)| \,. 
$$
Thus, $\gh$ is unimodular iff for all $x$ and $y$ in the same
orbit, $|S(x)| = |S(y)|$.
In particular, if $\gh$ is transitive, then $\gh$ is unimodular iff $|S(x)| =
|S(y)|$ for all neighbors $x$ and $y$.

Unimodularity of $\Aut(\gh)$ is a simple and
natural combinatorial property, as shown by \ref b.Schlichting/ and
\ref b.Trofimov/. Namely, if
$|\cbuldot|$ denotes cardinality
(for subsets of $\gh$) as well as Haar measure (for subsets of $\Aut(\gh)$),
then for any vertices $x, y \in \gh$,
$$
|S(x)y| / |S(y)x| = |S(x)| / |S(y)|\,; 
$$
thus, $\gh$ is unimodular iff for all $x$ and $y$ in the same orbit,
$$
|S(x)y| = |S(y)x| \,. \label e.unimod-schl-trof
$$
If $\gh$ is transitive, then $\gh$ is unimodular iff \ref
e.unimod-schl-trof/ holds for all neighbors $x, y$.

An {\bf end} of a graph $\gh$ is an equivalence class of infinite
nonself-intersecting paths in $\gh$, with two paths equivalent if for
all finite $A \subset \gh$, the paths are eventually in the same connected
component of $\gh \setminus A$.

\procl x.tend
Let $\gh$ be the regular tree of degree 3. Fix an end $\xi$
of $\gh$ and let $\gp$ be the set of automorphisms preserving $\xi$.
Then $\gp$ is a closed transitive subgroup of $\Aut(\gh)$
that is not unimodular.
For an example of a transitive graph $\gh$ whose full automorphism group
is not unimodular, add to the above tree, for each vertex $x$, the edge
between $x$ and its $\xi$-grandparent.  These examples were
described by \ref b.Trofimov/.
\endprocl

Next, we review amenability. Let $\gp$ be any locally compact group and
$L^\infty(\gp)$ be the Banach space of measurable real-valued functions on
$\gp$ that are essentially bounded with respect to Haar measure. A linear
functional on $L^\infty(\gp)$ is called a {\bf mean} if it maps the
constant function $\I{}$ to the number 1 and nonnegative functions to
nonnegative numbers. If $f \in L^\infty(\gp)$ and $\gpe \in \gp$, we write
$L_\gpe f(h):= f(\gpe h)$. We call a mean $\mu$ {\bf invariant} if
$\mu(L_\gpe f) = \mu(f)$ for all $f \in L^\infty(\gp)$ and $\gpe \in \gp$.
Finally, we say that $\gp$ is {\bf amenable} if there is an invariant mean
on $L^\infty(\gp)$.
F\o lner (see \ref b.Paterson/, Theorem 4.13) showed that $\gp$ is amenable
iff for every nonempty compact $B \subset \gp$ and $\epsilon > 0$, there is
a nonempty compact set $A \subset \gp$ such that $|B A \triangle A| \le
\epsilon |A|$. In this case, one often refers informally to $A$ as a F\o
lner set.

Now let $\gh = (\vertex, \edge)$ be a graph.
Given a set $K\subseteq \vertex$, let
$$
\mod K := \sum_{x\in K} |S(x)| \,.
$$
Note that $\mod \cbuldot$ is just counting measure if $\gh$ is unimodular
and Haar measure is normalized so that $|S(\bp)| = 1$.
Say that a
transitive graph $\gh$ is {\bf $\mod\cbuldot$-amenable} if for all $\epsilon
>0$, there is a finite $K\subset \vertex$ such that $\mod{\bdv K}<\epsilon
\mod K$.
If $\gh$ is unimodular, then this concept is the same as amenability of
$\gh$.
A mean on $\ell^\infty(\vertex)$ is called {\bf invariant} if
every $f \in \ell^\infty(\vertex)$ has the same mean as does $L_\gpe f$ (defined
as the function taking $x \mapsto f(\gpe x)$) for every $\gamma \in
\Aut(\gh)$.

For automorphism groups of graphs, amenability has the following
interpretations:

\procl t.eqvAmen  \procname{\ref B.BLPSgip/}
Let $\gh$ be a transitive graph. The following are equivalent:
\beginitems
\itemrm{(i)} $\Aut(\gh)$ is amenable;
\itemrm{(ii)} $\gh$ has an invariant mean;
\itemrm{(iii)} $\gh$ is $\mod\cbuldot$-amenable.
\enditems
\endprocl

\procl t.SW \procname{\ref B.SoardiWoess:rws/}
Let $\gh$ be a transitive graph.
Then $\gh$ is amenable iff $\Aut(\gh)$ is amenable and unimodular.
\endprocl

As usual, for any set $A$,
we write $2^A$ for $\{0, 1\}^A$ and identify it with the
collection of subsets of $A$.
It is given the usual product topology and Borel $\sigma$-field.

Let $\bp$ be a fixed vertex in $\gh$. In case $\gh$ is a tree, then $\bp$
will always designate the {\bf root} of $\gh$.
We denote by $|x|$ the graph distance
between $\bp$ and $x$ in $\gh$ for $x \in \vertex$. Let $B(x, n)$ denote
the set of vertices in $G$ within distance $n$ of $x$. Write
$$
\gr(\gh) := \liminf_{n \to\infty} |B(\bp, n)|^{1/n}
$$
for the (lower exponential) {\bf growth rate} of $\gh$.
When $\gh$ is transitive, we could replace $\liminf$ by $\lim$
because
$$
|B(\bp, m+n)| \le |B(\bp, m)| \cdot |B(\bp, n)|
\,.
$$
For any graph $\gh$, the fact that $\isov(\gh) \le |\bdv B(\bp,
n)|/|B(\bp, n)|$ for each $n$ implies that $1 + \isov(\gh) \le \gr(\gh)$.
In particular, if $\gh$ is nonamenable of bounded degree, then $\gr(\gh) > 1$.

We shall sometimes have processes indexed by elements of a graph as well as
by time. In order to distinguish between invariance under graph
automorphisms and under time, we shall reserve the term {\bf invariant} for
the former and use {\bf stationary} for the latter.

\bsection{Bernoulli Percolation on Transitive Graphs}{s.perc}

In {\bf Bernoulli($p$) bond percolation} on a graph,
each edge is {\bf open} (or {\bf occupied} or {\bf retained})
with probability $p$ independently.  Those edges that
are not open are {\bf closed} (or {\bf vacant} or {\bf removed}).
The corresponding product measure on $2^\edge$
is denoted $\P_p$.  The {\bf percolation subgraph} is the random
graph whose vertices are $\vertex$ and whose edges are the open edges.
Let $K(x)$ be the {\bf cluster} of $x$, that is, the
connected component of $x$ in the percolation subgraph.
We write
$$
\theta_x(p) :=\P_p\big[K(x)\hbox{ is infinite} ]\,.
$$
On a transitive graph, the value of $\theta_x(p)$ is independent of the
choice of $x$, whence the subscript $x$ is dropped.
The event that $K(x)$ is infinite is often written $x \leftrightarrow
\infty$. We also write $x \leftrightarrow y$ for $y \in K(x)$ and
$$
\tau_p(x, y) := \P_p[x \leftrightarrow y]
\,.
$$
Let
$$
p_c := p_c(\gh) := \inf\big\{p \st \theta(p)>0\big\}
$$
be the critical probability for percolation.

Bernoulli site percolation is
defined similarly with vertices replacing edges.
We shall use the superscripts ``bond" and ``site" when needed to
distinguish the two models.
See \ref b.Grimmett:newbook/
for more information about Bernoulli percolation.

If $\gh$ is a regular tree, then $K(\bp)$ is a Galton-Watson tree (except
for the first generation), so its analysis is easy and well known.
The first analysis of percolation on a nonamenable graph that is not a tree
was carried out by \ref b.GN:treeZ/ on the cartesian product of
the integers and a regular tree of sufficiently high degree.
They proved that for some
$p>p_c$, multiple infinite clusters coexist, while for other $p$, there is
a unique infinite cluster.
As a consequence of a method for studying random walks, \ref b.Lyons:groups/
gave a threshold for Bernoulli percolation on transitive graphs of exponential
growth (\ref t.LyGp/ below).
There followed the paper of \ref b.pyond/, which has spawned a considerable
amount of continuing research.

The results that follow are valid for both bond and site percolation when not
otherwise stated.
The only relations we shall state between site and bond percolation follow
from the usual coupling of the two processes (see, e.g., \ref
b.GrimStac:strict/ for the coupling): $\pcb(\gh) \le \pcs(\gh)$ for every
graph $\gh$, with strict inequality for most transitive $\gh$ proved by \ref
b.GrimStac:strict/; and $\pub(\gh) \le \pus(\gh)$ for transitive $\gh$,
where $p_u$ is defined in \ref e.pudef/ and \ref t.pu/(i) is being used to
establish the inequality.

It is well known that if $\gh$ is any infinite graph with the degree of each
vertex at most $d$, then $p_c(\gh) \ge 1/(d-1)$.
In the other direction, \ref b.Lyons:groups/ observed:

\procl t.LyGp
If $\gh$ is any transitive graph, then $p_c(\gh) \le 1/\gr(\gh)$.
\endprocl

This also follows immediately from

\procl t.sharp \procname{\ref B.AB:sharp/}
If $\gh$ is any transitive graph and $p < p_c(\gh)$, then $\E_p[|K(\bp)|] <
\infty$.
\endprocl

\noindent[\ref b.AB:sharp/ worked only on $\Z^d$, but their proof works in
greater generality.]

In particular, if $\gh$ is nonamenable and transitive, then it has exponential
growth, so that $0 < p_c(\gh) < 1$. The fact that $p_c(\gh) < 1$
was extended to nonamenable nontransitive graphs by \ref b.pyond/, who showed

\procl t.pciso
For any graph $\gh$, we have $\pcb(\gh) \le 1/\big(1+\isoe(\gh)\big)$ and
$\pcs(\gh) \le 1/\big(1+\isov(\gh)\big)$.
\endprocl

\procl r.anchored
The proof of \ref t.pciso/ actually gives better bounds, with $\isoe(\gh)$ and
$\isov(\gh)$ replaced by
$$
\ance(\gh) := \lim_{n \to \infty}\,
\inf \left\{{| \bde K| \over |K|} \st
o\in K\subset\vertex,\,K \hbox{ is connected},\,n\leq |K|<\infty \right\}
$$
and
$$
\ancv(\gh) := \lim_{n \to \infty}\,
\inf \left\{{| \bdv K| \over |K|} \st
o\in K\subset\vertex,\,K \hbox{ is connected},\,n\leq |K|<\infty \right\}
\,,
$$
respectively, the {\bf anchored expansion constants} introduced in \ref
b.BLS:pert/.
\endprocl

%
%
%

The next question concerns the number of infinite clusters when there is at
least one. When $\gh$ is
transitive, the argument of \ref b.NS/ shows that for any $p$, the number
of infinite clusters in Bernoulli($p$) percolation is an a.s.\ constant,
either 0, 1, or $\infty$.  As $p$ increases from 0 to 1, this constant goes
from 0 to $\infty$ to 1, possibly skipping $\infty$, as was shown by \ref
b.HP:mono/ in the unimodular case and by \ref b.Sch:sic/ in general:

\procl t.uniq
Let $\gh$ be a transitive graph.
Let $p_1 < p_2$.
If there is a unique infinite cluster $\P_{p_1}$-a.s., then
there is a unique infinite cluster $\P_{p_2}$-a.s.  Furthermore, in the
standard coupling of Bernoulli percolation processes, if there
exists an infinite cluster $\P_{p_1}$-a.s., then a.s.\ every
infinite $p_2$-cluster contains an infinite $p_1$-cluster.
\endprocl

Here, we refer to the {\bf standard coupling} of Bernoulli($p$) percolation
for all $p$ where, for bond percolation, say, each edge $e \in \edges$ is
assigned an independent uniform $[0, 1]$ random variable $U(e)$ and the
edges where $U(e) \le p$ are retained for Bernoulli($p$) percolation.

If we define
$$
p_u(\gh) := \inf \big\{ p \st \hbox{there is a unique infinite cluster in
Bernoulli($p$) percolation} \big\}
\,,
\label e.pudef
$$
then it follows from \ref t.uniq/ that when $\gh$ is transitive,
$$
p_u(\gh) = \sup \{ p \st \hbox{there is not a unique infinite cluster in
Bernoulli($p$) percolation} \}\,.
$$

It is not hard to show that when $\gh$ is a transitive graph with at least 3
ends, then $p_u(\gh) = 1$.
Since nonamenable transitive graphs cannot have only two ends, the
remaining cases fall under the following conjecture, suggested in a
question of \ref b.pyond/:

\procl g.punontriv
If $\gh$ is a transitive nonamenable graph with one end, then $p_u(\gh) < 1$.
\endprocl

This conjecture has been confirmed in the following cases:
\beginbullets

$\gh$ is a Cayley graph of a finitely presented group (\ref B.BabB:1end/);

$\gh$ is planar (\ref b.Lalley:Fuch/ for site percolation on co-compact
Fuchsian groups of genus at least 2, and \ref b.BS:hp/ for percolation
in general; the full result can also be
deduced from the argument of \ref b.BabB:1end/);

$\gh$ is the cartesian product of two infinite graphs (\ref B.HPS:merge/);

$\gh$ is a Cayley graph of a Kazhdan group, i.e., a group with Kazhdan's
property T (\ref B.LS:indis/).
\endbullets

Some additional information about the uniqueness phase is contained in the
following theorem:

\procl t.pu
Let $\gh$ be a transitive graph.
\beginitems
\itemrm{(i)}
\procname{\ref B.Sch:sic/}
$$
p_u(\gh) = \inf \{ p \st \sup_R \inf_x \P_p[B(o, R) \leftrightarrow B(x, R)]
= 1\}\,.
\label e.puBB
$$
\itemrm{(ii)}
\procname{\ref B.LS:indis/}
If $\gh$ is unimodular and $\inf_x \tau_p(o, x) > 0$,
then there is a unique infinite cluster $\P_p$-a.s. Therefore,
$$
p_u(\gh) = \inf \{ p\st \inf_x \tau_p(o, x) > 0 \}
\,.
\label e.putau
$$
\enditems
\endprocl

Equation \ref e.putau/ implies \ref e.puBB/, but it is unknown whether \ref
e.putau/ holds in the nonunimodular case.

As is well known, when $\gh$ is amenable and transitive, there can never
be infinitely many infinite clusters (\ref b.BK:uni/ for $\Z^d$ and
\ref b.GKN:uni/ in general),
whence $p_c(\gh) = p_u(\gh)$.
Behavior that is truly different from the amenable case arises when there
{\it are\/} infinitely many infinite clusters. This has been conjectured
always to be the case on nonamenable transitive graphs for an interval of
$p$:

\procl g.pcpu \procname{\ref B.pyond/}
If $\gh$ is a transitive nonamenable graph, then $p_c(\gh) < p_u(\gh)$.
\endprocl

This has been confirmed in certain cases:
\beginbullets

if $\gh$ is the product of any transitive graph
with a regular tree of sufficiently high degree
(\ref b.GN:treeZ/ when the transitive graph is $\Z$ and \ref b.pyond/ in
general);

if $\gh$ is planar (\ref b.Lalley:Fuch/ for site percolation on co-compact
Fuchsian groups of genus at least 2, and \ref b.BS:hp/ for percolation
in general);

for bond percolation if $\isoe(\gh)/\gdeg \ge 1/\sqrt 2$ and for site
percolation if $\isov(\gh)/\gdeg \ge 1/\sqrt 2$, where $\gdeg$ is
the degree of $\gh$ (\ref B.RS:isoD/); this implies the first bulleted case
above.

if $\gh$ is any Cayley graph of a group of cost larger than 1. This
includes, first,
free groups of rank at least 2 and
fundamental groups of compact surfaces of genus larger than 1.
Second, let $\gp_1$ and $\gp_2$ be two groups of finite cost
with $\gp_1$ having cost larger than 1.
Then every amalgamation of $\gp_1$ and $\gp_2$ over an amenable group has
cost larger than 1. Third, every
HNN extension of $\gp_1$ over an amenable group has cost larger than 1.
For the definition of cost and proofs that these groups have cost larger
than 1, see \refbmulti {Gaboriau:mercuriale,Gaboriau:cost}. The proof that
$p_c(\gh) < p_u(\gh)$ follows fairly easily from \ref t.death/ below.
\endbullets

The third bulleted case above uses the following lower bound for $p_u(\gh)$
[or the weaker bound of \ref b.pyond/, Theorem 4]. Here,
a {\bf simple cycle} is a cycle that does not use any vertex or edge
more than once.

\procl t.pulower \procname{\ref B.Schramm:pc97/}
Let $\gh$ be a transitive graph and let $a_n(\gh)$ be the number of simple
cycles of length $n$ in $\gh$ that contain $\bp$. Then
$$
p_u(\gh) \ge \liminf_{n \to\infty} a_n(\gh)^{-1/n}
\,.
\label e.pulower
$$
\endprocl

\proof
We give the proof for site percolation, the proof for bond percolation
being similar.
Let $U(x)$ be independent uniform $[0, 1]$ random variables indexed by
$\vertex$.
Take $p > p' > p_u \ge p_c$.  In order to show that $p_u(\gh) \ge \liminf_{n
\to\infty} a_n(\gh)^{-1/n}$, we shall show that $\sum_n a_n(\gh) p^n =
\infty$.
Let $\omega$ be the open subgraph formed by the vertices $x$
with $U(x) \le p$.
First, observe that since $\omega$ contains a.s.\ a unique infinite
cluster, that infinite cluster $K$ has only one end, since otherwise
removing a finite number of edges would create more than one infinite
cluster.

Second, with positive probability, there are two (edge- and vertex-)
disjoint infinite rays in $K$.  Otherwise, by Menger's theorem, for any
vertex $x \in K$, a.s.\ there would be infinitely many vertices $x_n$, each
of whose removal would leave $x$ in a finite open component. But given
$\omega$, given any such vertex $x$, and given any such vertices $x_n$,
$U(x_n) > p'$ a.s.\ for infinitely many $n$.  This means that $K(x)$ is
finite $\P_{p'}$-a.s.\ and contradicts $p'>p_c$.

Therefore, with positive probability there are two infinite rays in
$\omega$ starting
at $\bp$ that are disjoint except at $\bp$.
Since $K$ has only one end, the two
rays may be connected by paths in $\omega$ that stay
outside arbitrarily large balls.
In particular, there are an infinite number of simple cycles in $\omega$
through $\bp$, whence the expected number of simple cycles through $\bp$ in
$\omega$ must be infinite.
That is, $\sum_n a_n(\gh) p^n = \infty$.
\Qed

Additional evidence for \ref g.pcpu/ is provided by

\procl t.PSN \procname{\ref B.PakSN:uniq/}
For any finitely generated nonamenable group $\gp$, there exists some
Cayley graph $\gh$ of \/ $\gp$ with $p_c(\gh) < p_u(\gh)$.
\endprocl

\noindent The proof of \ref t.PSN/ shows that the Cayley graph can be found
so as to satisfy Schonmann's condition above that $\iso(\gh)/d \ge 1/\sqrt
2$.

We next discuss behavior of percolation at the critical values $p_c$ and
$p_u$.
It has been long conjectured that there are no infinite clusters at the
critical value $p_c(\gh)$ when $\gh$ is a euclidean lattice, i.e.,
$\theta\big(p_c(\gh)\big) = 0$. \ref b.pyond/ extended this conjecture to all
transitive $\gh$. This was confirmed in the unimodular nonamenable case:

\procl t.death \procname{\refbmulti{BLPSgip,BLPSdeath}}
If $\gh$ is a unimodular nonamenable transitive graph, then
$\theta\big(p_c(\gh)\big) = 0$.
\endprocl

It follows from Theorems \briefref t.uniq/, \briefref t.death/, and a result
of \ref b.BK:conti/ that $\theta(p)$ is continuous in $p$ on each
nonamenable unimodular transitive $\gh$.
(This was proved earlier by \ref b.Wu:CtnyHyp/ for a graph that is not
transitive but is similar to the hyperbolic plane.)

It is unknown how many infinite clusters there are at $p_u$.  It is known
that there is a unique infinite cluster at $p_u$ when $\gh$ is planar,
nonamenable and transitive (\ref B.BS:hp/). On the other hand, there cannot
be exactly one infinite cluster at $p_u$ when $\gh$ is a cartesian product
(of infinite transitive graphs) with a nonamenable automorphism group (\ref
b.S:nuni/ in the case of a tree cross $\Z$ and \ref b.Peres:pu/ in general)
or when $\gh$ is a Cayley graph of a Kazhdan group (due to Peres; see \ref
b.LS:indis/).

Finally, we discuss briefly the nature of the infinite clusters when there
are infinitely many of them; see \ref b.BLS:pert/ and \ref b.HSS:dilute/
for more on this topic. A basic
result is that when there are infinitely many infinite clusters, they are
``indistinguishable" from each other:

\procl t.indist
Let $\gh$ be a transitive unimodular graph.
Let $\ev A$ be a Borel measurable set of subgraphs of $\gh$
that is invariant under the automorphism group of $\gh$.
Then either $\P_p$-a.s.\ all infinite clusters are in $\ev A$,
or $\P_p$-a.s.\ they are all outside of $\ev A$.
\endprocl

\noindent \ref t.indist/ was proved by \ref b.HP:mono/ for increasing
sets $\ev A$ and for all but possibly one value of $p$, while it was
proved in general (and for certain other percolation processes) by \ref
b.LS:indis/.

For example, $\ev A$ might be the collection of all transient
subgraphs of $\gh$, or the collection of all subgraphs that have
a given asymptotic rate of growth, or the collection of all subgraphs
that have no vertex of degree $5$.

If $\ev A$ is the collection of all transient subgraphs
of $\gh$, then \ref t.indist/ shows that almost surely, either all infinite
clusters of $\omega$ are
transient [meaning that simple random walk on them is transient],
or all clusters are recurrent.
In fact, \ref b.LS:indis/ show that if $\gh$ is nonamenable,
then a.s.\ all infinite clusters are transient if Bernoulli percolation
produces more than one infinite component.  (\ref b.BLS:pert/ show
that the same is true if Bernoulli percolation produces a single
infinite component.)


We illustrate some uses of \ref t.indist/ by proving two theorems (though
it should be noted that the original direct proofs of these theorems are
simpler than the proof of \ref t.indist/). \ref t.indist/ is also used to
prove \ref t.pu/(ii).

\smallbreak\noindent{\it Proof of \ref t.uniq/ in the unimodular
case.}\enspace
Suppose that there exists an infinite cluster $\P_{p_1}$-a.s.
Let $\omega$ be the open subgraph of the $\P_{p_2}$ process and let $\eta$ be
an independent $\P_{p_1/p_2}$ process.
Thus, $\omega \cap \eta$ has the law of $\P_{p_1}$ and, in fact, $(\omega \cap
\eta, \omega)$ has the same law as the standard coupling of $\P_{p_1}$ and
$\P_{p_2}$.
By assumption, $\omega \cap \eta$ has an infinite cluster a.s.
Thus, for some cluster $C$ of $\omega$, we have $C \cap \eta$ is infinite with
positive probability, hence, by Kolmogorov's 0-1 law, with probability 1.
By \ref t.indist/, this holds for every cluster $C$ of $\omega$.
\Qed

An extension to \ref t.uniq/ in the unimodular case is as follows. It is
unknown whether it holds in the nonunimodular case.

\procl t.merge \procname{\ref B.HPS:merge/}
Let $\gh$ be a transitive unimodular graph.
Let $p_1 < p_2$ be such that there are infinitely many infinite clusters
$\P_{p_1}$-a.s.\ and $\P_{p_2}$-a.s. In the standard coupling of
Bernoulli percolation processes on $\gh$,
a.s.\ every infinite $p_2$-cluster contains infinitely many infinite
$p_1$-clusters.
\endprocl

\proof (due to R.~Schonmann)
The number of infinite $p_1$-clusters contained in a $p_2$-cluster is a
random variable whose distribution is the same for each infinite
$p_2$-cluster by \ref t.indist/. In fact, by an extension of \ref t.indist/
involving random scenery that is stated by \ref b.LS:indis/ (and that has
the same proof), this random variable is constant a.s. Thus, each infinite
$p_2$-cluster has the same number of infinite $p_1$-clusters a.s. Since two
infinite $p_2$-clusters could merge through the addition of finitely many
edges, the number of infinite $p_1$-clusters contained in an infinite
$p_2$-cluster could change unless that number were infinite. \Qed

\ref t.indist/ does not hold for nonunimodular graphs (\ref B.LS:indis/).
However, \ref b.HPS:merge/ have found a replacement
that does hold without the unimodularity assumption (as long as $p>p_c$;
presumably, this caveat is not important since presumably there are no
infinite clusters at $p_c$).
Define $\ev A$ to be {\bf robust} if for every infinite connected
subgraph $C$ of $\gh$ and every edge $e \in C$, we have
$C \in \ev A$ iff there is an infinite connected component of
$C \setminus\{e\}$ that lies in $\ev A$.
For example, transience is a robust property.

\procl t.robust \procname{\ref B.HPS:merge/}
Let $\gh$ be a transitive graph.
Let $p > p_c(\gh)$
and let $\ev A$ be a robust Borel measurable set of subgraphs of $\gh$.
Assume that $\ev A$ is invariant under the automorphism group of $\gh$.
Then either $\P_p$-a.s.\ all infinite percolation components are in $\ev A$,
or $\P_p$-a.s.\ they are all outside of $\ev A$.
\endprocl

Finally, it should be noted that \ref b.pyond/ contains several interesting
questions and conjectures about various families of graphs, including
nontransitive and amenable graphs. One may consult \ref b.pyond:recent/ for
updates concerning progress on Bernoulli percolation on general graphs.

\bsection{Invariant Percolation on Transitive Graphs}{s.invar}

As we have mentioned in the introduction, there are interesting and useful
results about invariant percolation, especially on transitive nonamenable
graphs. We give a sample of these results here that show their nature and
how they can be used. In addition, we illustrate the powerful
mass-transport technique.
All formally stated results in this section are from \BLPSgip, which will
be referred to simply as \BLPSgip\ throughout this section.

A {\bf bond percolation process} is a pair $(\P,\omega)$,
where $\omega$ is a random element in $2^\edges$ and $\P$ denotes
the distribution (law) of $\omega$.
We shall say that $\omega$ is the {\bf configuration} of the percolation.
A {\bf site percolation process} $(\P,\omega)$ is given by a probability
measure $\P$ on $2^{\vertex(G)}$, while a (mixed) {\bf percolation}
is given by a probability measure on $2^{\vertex(G)\cup\edges(G)}$
that is supported on subgraphs of $G$.
If $\omega$ is a bond percolation process,
then $\hat\omega :=\vertex(G)\cup\omega$
is the associated mixed percolation.
In this case, we shall not distinguish between $\omega$ and
$\hat\omega$, and think of $\omega$ as a subgraph of $G$.
Similarly, if $\omega$ is a site percolation, there is an
associated mixed percolation
$\hat\omega:=\omega\cup\bigl(\edges(G)\cap(\omega\times\omega)\bigr)$,
and we shall not bother to distinguish between $\omega$ and
$\hat\omega$.

If $x\in\vertex(G)$ and $\omega$ is a percolation on $G$,
the {\bf cluster} (or {\bf component}) $K(x)$ of $x$ in $\omega$
is the set of vertices in $\vertex(G)$ that can be connected
to $x$ by paths contained in $\omega$.
We shall not distinguish between the cluster $K(x)$ and
the graph
$\Bigl(K(x),\bigl(K(x)\times K(x)\bigr)\cap\omega\Bigl)$ whose vertices
are $K(x)$ and whose edges are the edges in $\omega$ with endpoints
in $K(x)$.

A percolation process $(\P,\omega)$ in a graph $\gh$ is called {\bf
invariant} if $\P$ is invariant under $\Aut(\gh)$.  Invariant percolation
has proved useful for the study of Bernoulli percolation as well as
other processes such as the random cluster model, as we shall see below. It
is also interesting in itself.

We first present the very useful Mass-Transport Principle.
Early forms of the mass-transport method were used by \ref b.Adams:trees/
and \ref b.BergMees:stab/.
It was introduced in the study of percolation by \ref b.Hag:deptree/ and
developed further in \BLPSgip.
Let $\xi$ be some (automorphism-)invariant process on $\gh$, such as
invariant percolation, and let $F(x, y;\xi)\in[0,\infty]$ be a function of
$x, y \in\vertex$ and $\xi$.  Suppose that $F$ is invariant under the
diagonal action of $\Aut(\gh)$; that is, $F( \gpe x,  \gpe y; \gpe
\xi)=F(x, y,\xi)$ for all $ \gpe  \in \Aut(\gh)$.  We think of giving
each vertex $x \in \vertex$ some initial mass, possibly depending on
$\xi$, then redistributing it so that $x$ sends $y$ the mass $F(x, y;
\xi)$. With this terminology, one hopes for ``conservation'' of mass, at
least in expectation.  Of course, the total amount of mass is usually
infinite. Nevertheless, there is a sense in which mass is conserved; in the
transitive unimodular setting, we have that the expected mass at a vertex
before transport equals the expected mass at a vertex afterwards.  More
generally, mass needs to be weighted according to the Haar measure of the
stabilizer. Since $F$ enters only in expectation, it is
convenient to set $f(x, y) := \E F(x, y;\xi)$.
For the reader to whom this is new, it is recommended to consider only the
unimodular case; then all factors of $|S(x)|$ become 1 and all $*$'s below
can be omitted.

\proclaim Mass-Transport Principle. 
If $\gh$ is a transitive graph and
$f : \gh \times \gh \to [0, \infty]$ is invariant
under the diagonal action of $\Aut(\gh)$, then
$$
\sum_{x \in \vertex} f(\bp, x)
= \sum_{x \in \vertex} f(x, \bp) |S(x)|/|S(\bp)| \,.
$$

For a subgraph $K \subset \gh$, let $\deg_K(x)$ denote the degree of
$x$ in $K$. If $K$ is finite and nonempty, put
$$
\alpha_K^* := {1 \over  \mod K} \sum_{x \in K} \deg_K(x) |S(x)|
\,;
$$
this is the average (internal) degree in $K$, appropriately weighted if the
graph is not unimodular.
Then define
$$
\alpha^*(\gh) := \sup \{ \alpha_K^* \st K \subset \gh \hbox{ is finite and
nonempty} \}\,.
$$
If $\gh$ is a regular graph of degree $\gdeg$, then
$$
\alpha^*(\gh) + \isoe^*(\gh) = \gdeg\,,
\label e.sumd
$$
where
$$
\isoe^*(\gh) := \inf \left\{ {1 \over \mod K} \sum_{[x,y] \in \bde K} |S(x)|
    \st K\subset \vertex \hbox{ is finite and nonempty} \right\} \,.
$$
For a random subgraph $\omega$ of $\gh$
and a vertex $x \in \gh$, define
$$
D^*(x) := \sum_{[x, y] \in \omega} |S(y)|/|S(x)|
\,.
$$
Let
$$
\gdeg^* := \sum_{[\bp, y] \in \gh} |S(y)|/|S(\bp)|
\,.
$$

We give two simple but useful applications of the Mass-Transport Principle
to illustrate the method. The first is quantitative, while the second is
qualitative. Both were proved earlier by \ref b.Hag:deptree/ for regular
trees. (His paper was the original impetus for \BLPSgip.)
Write
$$
\thetaP := \P[\bp \leftrightarrow \infty]
\,.
\label e.thetadef
$$

\procl t.NonuniTh 
Let $\gh$ be a nonamenable transitive graph and
$\P$ be an invariant bond percolation on $\gh$. Then
$$
\thetaP
\ge [\E D^*(\bp) - \alpha^*(\gh)]/\isoe^*(\gh)\,.
\label e.bound
$$
In particular, if $\E D^*(\bp) > \alpha^*(\gh)$, then
$\thetaP > 0$.
\endprocl

The intuition is that if the expected (weighted) degree of a vertex is
larger than the average internal degree of finite subgraphs, then it must
be carried by some infinite components.

\proof
Let $I_x$ be the indicator that $K(x)$ is finite.
We put mass $D^*(x)I_x$ at each $x \in \vertex$.
In each finite component, the masses are redistributed proportionally to
the weights $|S(y)|$ (for $y$ in the component) among the vertices in that
component. Since $\P$ is invariant, so is this mass transport. Formally, we
use the function
$$
f(x, y) := \E\left[ I_x \II{y \in K(x)} {D^*(x) |S(y)|\over \mod{K(x)}} \right]
\,,
$$
which is automorphism invariant.
We have
$$
\sum_{z \in \vertex} f(\bp, z)
= \E[D^*(\bp) I_\bp]\,.
$$
On the other hand,
\begineqalno
\sum_{y \in \vertex} f(y, \bp) |S(y)|/|S(\bp)|
&=
\Ebigg{I_\bp \!\! \sum_{y \in K(\bp)} {D^*(y) |S(\bp)|\over \mod{K(\bp)}}
{|S(y)| \over |S(\bp)|}}
\cr&= \Eleft{ \alpha^*_{K(\bp)} I_\bp }
\le  \alpha^*(\gh) \big(1 - \thetaP\big)\,.
\cr
\endeqalno
Since $D^* \le d$ everywhere, the Mass-Transport Principle implies that
\begineqalno
\E D^*(\bp) - \gdeg \thetaP
&\le \E[D^*(\bp) I_\bp]
=
\sum_{z \in \vertex} f(\bp, z)
=
\sum_{y \in \vertex} f(y, \bp) |S(y)|/|S(\bp)|
\cr&\le \alpha^*(\gh) \big(1 - \thetaP\big)\,.
\cr
\endeqalno
A little algebra using \ref e.sumd/ completes the proof. \Qed

Variations on this result have proved useful. For example \BLPSgip, if in
addition to the above hypotheses, $\gh$ is unimodular, $\P$ has the
property that all components are trees a.s., and $\E D^*(\bp) \ge 2$, then
$\thetaP > 0$.

Our second application of the Mass-Transport Principle helps us to
count the ends of the components in the configuration of a percolation that
is invariant under a unimodular automorphism group:

\procl p.nmbr-ends 
Let $\gh$ be a unimodular transitive
graph. Let $\omega$ be the configuration of an invariant percolation on
$\gh$ such that $\omega$ has infinite components with positive probability.
Almost surely every component of $\omega$ with at least 3 ends has
infinitely many ends.
\endprocl

\proof
Let $\omega_1$ be the union of the components $K$ of $\omega$
whose number $n(K)$ of ends is finite and at least $3$.
Given a component $K$ of $\omega_1$,
there is a connected subgraph $A\subset K$ with minimal $|\vertex(A)|$
such that $K\setminus A$ has $n(K)$ infinite components.
Let $H(K)$ be the union of all such subgraphs $A$.
It is easy to verify that any two such subgraphs $A$ must intersect,
and therefore $H(K)$ is finite.
Let $H(\omega_1)$ be the union of all $H(K)$, where $K$ ranges over the
components of $\omega_1$.

Begin with unit mass at each vertex $x$ that belongs to
a component $K$ of $\omega_1$, and transport it equally to the
vertices in $H(K)$.
Then the vertices in $H(\omega_1)$ receive infinite mass.
By the Mass-Transport Principle, no vertex can receive infinite
mass, which means that $\omega_1$ is empty a.s.  \Qed

Among the characterizations of amenability via invariant percolation that
appear in \BLPSgip, we single out one that relates to the absence of phase
transition:

\procl t.treechar 
Let $\gh$ be a
transitive graph. Then each of the following
conditions implies the next one:
\beginitems
\itemrm{(i)} $\gh$ is amenable;
\itemrm{(ii)} there is an invariant random nonempty subtree of $\gh$ with at
most 2 ends a.s.;
\itemrm{(iii)} there is an invariant random nonempty
connected subgraph $\omega$ of $\gh$ that satisfies $p_c(\omega) = 1$
 with positive probability;
\itemrm{(iv)} $\Aut(\gh)$ is amenable.
\enditems
If $\gh$ is assumed to be unimodular, then all four
conditions are equivalent.
\endprocl


To see one use of \ref t.treechar/, we present the proof of part of \ref
t.death/. (In fact, here we do not need the assumption of unimodularity.)

\procl c.Pc-Pu 
If $\gh$ is a transitive graph with a nonamenable
automorphism group and Bernoulli($p$) percolation
produces a unique infinite cluster a.s., then $p > p_c(\gh)$.
\endprocl

\proof Suppose that $p = p_c(\gh)$ and that there is a unique infinite
cluster a.s. Then the infinite cluster $K$ has
$p_c(K)$ = 1 a.s. Hence $\Aut(\gh)$ is amenable. \Qed

Next, we present a characterization of unimodularity in terms of the
expected degree of vertices in infinite components. Since any connected
finite graph with vertex set $\vertex$ has average degree at least
$2-2/|\vertex|$, one might expect that for invariant percolation on a
transitive graph $\gh$ with all components infinite a.s., the expected degree
of a vertex is at least 2. This inequality is true when $\gh$ is
unimodular, but surprisingly, whenever $\gh$ is not unimodular, there is an
invariant percolation where the inequality fails.

\procl t.modular 
Let $\gh$ be a transitive graph.  Let $m$ be
the minimum of $|S(x)|/|S(y)|$ for $x$, $y$ neighbors in $\gh$. Then for any
invariant percolation that yields infinite components with positive
probability, the expected degree of $\bp$ given that $\bp$
is in an infinite component is at least $1 + m$. This is sharp for
all $\gh$ in the sense that there is an invariant bond percolation
on $\gh$ with every vertex belonging to an infinite component and having
expected degree $1+m$.
\endprocl

A {\bf forest} is a graph all of whose components are trees. The following
theorem concerning phase transition on percolation components was shown
when $\gh$ is a tree by \ref b.Hag:deptree/.

\procl t.ends 
Let $\gh$ be a unimodular
transitive graph. Let $\omega$ be the configuration of
an invariant percolation on $\gh$
such that $\omega$ has infinite components with positive probability.  If
\beginitems
\itemrm{(i)} some component of $\omega$ has at least 3 ends with positive
probability,
\enditems
then
\beginitems
\itemrm{(ii)} some component of $\omega$ has $p_c < 1$ with positive
probability and
\itemrm{(iii)} $\Eleft{D^*(\bp) \bigm| |K(\bp)| = \infty} > 2$.
\enditems
If $\omega$ is a forest a.s., then the three conditions are equivalent.
\endprocl

To show how \ref t.ends/ can be used, we now complete the proof of \ref
t.death/. (A more direct proof of \ref t.death/ is provided by \ref
b.BLPSdeath/.)

\proofof t.death
Let $\omega$ be the configuration of critical Bernoulli percolation on $\gh$.
Then every infinite cluster $K$ of $\omega$ has $p_c(K) = 1$ a.s.
As we have mentioned, the number of infinite clusters of $\omega$ is
equal a.s.\ to 0, 1 or $\infty$.
\ref c.Pc-Pu/ rules out a unique infinite cluster.
If there were more than one infinite cluster, then
by opening the edges in a large ball, we
see that there would be, with positive probability, a cluster
with at least 3 ends. In light of  \ref t.ends/, this would mean
that with positive probability, some infinite cluster $K$
had $p_c(K)<1$.  This is a contradiction.
\Qed

\bsection{Ising, Potts, and Random Cluster Models on Transitive Graphs}{s.ising}

Ising and Potts models on graphs are defined using interaction strengths
along bonds (here assumed identically 1), Boltzmann's constant $k_B$, and
the temperature $T$.
These last two quantities always appear together in the expression $\beta
:= 1/(k_B T)$, called the {\bf inverse temperature}.
Given a finite graph $\gh$ and an integer $q \ge 2$, let $\omega \in \{1,
2, \ldots, q\}^\vertex$. Write $I_\omega(e)$ for the indicator that
$\omega$ takes different values at the endpoints of the edge $e$. The {\bf
energy} (or {\bf Hamiltonian}) of $\omega$ is
$$
H(\omega) :=
2\sum_{e \in \edge} I_\omega(e)
\,.
$$
The {\bf Potts measure} $\fpt(\beta) = \fpt^\gh(\beta)$ is the probability
measure on $\{1, 2, \ldots, q\}^\vertex$ that is proportional to
$e^{-\beta H(\omega)}$.
In the case $q=2$, it is more customary to use $\{-1, 1\}^\vertex$ in place
of $\{1, 2\}^\vertex$, and the measure is called the {\bf Ising measure}.

To define such measures on infinite graphs $\gh$, one can proceed via {\bf
exhaustions} of $\gh$, i.e., sequences of finite subgraphs $\gh_n$ that are
increasing and whose union is all of $\gh$. There are several ways to do
this, in fact, and crucial questions are whether some of the limits they
give are the same. One way to take a limit is simply to define
$\fpt^\gh(\beta)$ to be the weak${}^*$ limit of $\fpt^{\gh_n}(\beta)$; this
is called the {\bf free Potts measure} on $\gh$. Another way is as follows.
Let $\ppt_k^{\gh_n}(\beta)$ be the probability measure
$\fpt^{\gh_n}(\beta)$ conditioned on having $\omega(x) = k$ for every $x
\in \bdiv\gh_n$, where
$$
\bdiv K := \{ x \in K \st \texists{y \notin K} x \sim y\}
$$
denotes the {\bf internal vertex boundary} of any subset $K
\subset \vertex$.
Then define the {\bf Potts measure} $\ppt_k^\gh(\beta)$ to
be the weak${}^*$ limit of $\ppt_k^{\gh_n}(\beta)$.
These limits always exist (see, e.g., \ACCN, referred to later as \ACCN).
It will be convenient to define the {\bf wired Potts measure}
$\wpt^\gh(\beta)$ to be $\sum_{k=1}^q \ppt_k^\gh(\beta)/q$.
Note that if $\gh_n^*$ denotes the graph obtained from $\gh_n$ by identifying
all of the vertices in $\bdiv \gh_n$ to a single vertex, then
$\wpt^\gh(\beta)$ is the weak${}^*$ limit of $\fpt^{\gh_n^*}(\beta)$.

To define Potts measures in general, write $\omega \restrict \vertex'$ for
the restriction of $\omega$ to $\vertex' \subset \vertex$.
For a finite subset $\vertex' \subset \vertex$, let $\gh'$ denote
the subgraph of $\gh$ induced by $\vertex'$, i.e., $\gh' := \big(\vertex',
(\vertex' \times \vertex') \cap \edge\big)$.
For $\omega' \in \{1, \ldots, q\}^{\vertex'}$, write $\bd \omega' :=
\omega' \restrict \bdiv \gh'$.
We call $\P$ a {\bf Potts measure} on $\gh$ at inverse temperature $\beta$
if $\P$ is a Markov random field and for all finite $\vertex' \subset
\vertex$ and all $\omega' \in \{1, \ldots, q\}^{\vertex'}$,
$$
\PBig{\omega \restrict \vertex' = \omega' \Bigm| \omega \restrict
\bdiv \gh' = \bd \omega'}
=
\fpt^{\gh'}(\beta)\Big[\omega = \omega' \Bigm| \omega
\restrict \bdiv \gh' = \bd \omega'\Big] \,.
$$
It is easy to verify that the measures $\fpt^\gh(\beta)$ and
$\ppt_k^\gh(\beta)$ are Potts measures in this sense.

Potts measures are intimately connected to random cluster measures,
introduced by \ref b.FK/ and \refbmulti{Fortuin2,Fortuin3}.  See
\ref b.OH:RCphase/ for a survey of the relationships and \ref b.Grimmett:RC/
for more details on random cluster measures, especially on $\Z^d$.
Random cluster measures depend on two parameters, $p \in (0, 1)$ and $q > 0$.
We restrict ourselves to $q \ge 1$ since the measures with $q < 1$ behave
rather differently and are poorly understood; they are also unrelated to
Potts measures.
Given a finite graph $\gh$ and $\omega \in 2^\edge$, write $\ncomp(\omega)$
for the number of components of $\omega$.
The {\bf random cluster measure} with parameters $(p, q)$ on $\gh$, denoted
$\frc(p, q) = \frc^\gh(p, q)$, is the probability measure on $\edge$
proportional to $q^{\ncomp(\omega)} \P_p(\omega)$, i.e., the Bernoulli($p$)
percolation measure $\P_p$ biased by $q^{\ncomp(\omega)}$ (and
renormalized).
On infinite graphs $\gh$, there are again several ways to define random
cluster measures. The ones that concern us are obtained by taking limits
over exhaustions $\gh_n$ of $\gh$. Namely, define $\frc^\gh(p, q)$ to
be the weak${}^*$ limit of $\frc^{\gh_n}(p, q)$; this is called the {\bf
free random cluster measure} on $\gh$.
Define the {\bf wired random cluster measure} $\wrc^\gh(p, q)$ to be
the weak${}^*$ limit of $\frc^{\gh_n^*}(p, q)$.
These limits always exist (see, e.g., \ACCN).
Furthermore, they have positive correlations and so the free random cluster
measure is stochastically dominated by the wired random cluster measure \ACCN.

Note that there is another use of ``wired'' in the literature, although
when $\gh$ has only one end, the meaning is the same as the present one.
In the terminology of \ref b.Grimmett:RC/, the above random cluster
measures are ``limit random cluster measures''.  We do not examine whether
they satisfy so-called Gibbs specifications.

Since all the above limits exist regardless of the exhaustion chosen, the
limiting measures are invariant under all graph automorphisms.

The fundamental relation between Potts and random cluster measures is the
following: Let $\gh$ be a finite
graph and $q \ge 2$ an integer. Suppose that $p = 1-e^{-2\beta}$.
\item{$\bullet$}
If $\omega \in \{1, \ldots, q\}^\vertex$ is chosen with distribution
$\fpt(\beta)$ and $\eta \in 2^\edge$ is chosen independently
with distribution $\P_p$,
then $(1 - I_\omega) \eta$ has the distribution $\frc(p, q)$.
\item{$\bullet$}
Choose $\eta \in 2^\edge$ with distribution $\frc(p, q)$. For each
component of $\eta$, choose independently and uniformly an element of
$\{1, \ldots, q\}$, assigning this element to every vertex in that
component. The resulting $\omega \in \{1, \ldots, q\}^\vertex$ has
distribution $\fpt(\beta)$.

\noindent See \ACCN\ or \ref b.OH:RCphase/ for proofs.
By taking weak${}^*$ limits and using positive correlations, one obtains
corresponding statements for infinite graphs (see the proof of Theorem
2.3(c) of \ACCN): First, the two above statements hold as written for
infinite graphs.
Second:
\item{$\bullet$}
If $\omega \in \{1, \ldots, q\}^\vertex$ is chosen with distribution
$\wpt(\beta)$ and $\eta \in 2^\edge$ is chosen independently
with distribution $\P_p$,
then $(1 - I_\omega) \eta$ has the distribution $\wrc(p, q)$.
\item{$\bullet$}
Choose $\eta \in 2^\edge$ with distribution $\wrc(p, q)$. For each
component of $\eta$, choose independently and uniformly an element of
$\{1, \ldots, q\}$, assigning this element to every vertex in that
component, where all infinite components are regarded as a single component
(``connected at infinity").
The resulting $\omega \in \{1, \ldots, q\}^\vertex$ has
distribution $\wpt(\beta)$.
\smallbreak\noindent
Third:
\item{$\bullet$}
If $\omega \in \{1, \ldots, q\}^\vertex$ is chosen with distribution
$\ppt_k(\beta)$ and $\eta \in 2^\edge$ is chosen independently
with distribution $\P_p$,
then $(1 - I_\omega) \eta$ has the distribution $\wrc(p, q)$.
\item{$\bullet$}
Choose $\eta \in 2^\edge$ with distribution $\wrc(p, q)$. For each
finite component of $\eta$, choose independently and uniformly an element of
$\{1, \ldots, q\}$, assigning this element to every vertex in that
component.
Assign each vertex in an infinite component the color $k$.
The resulting $\omega \in \{1, \ldots, q\}^\vertex$ has
distribution $\ppt_k(\beta)$.

\smallbreak

Recall the notation \ref e.thetadef/. From the preceding relations, we obtain:

\procl p.PtRC
Let $\gh$ be any graph and $q \ge 2$ an integer. Let $\beta > 0$ and $p := 1
- e^{-2\beta}$. Then
\beginitems
\itemrm{(i)} \procname{\ref B.Jonasson:RCamen/}
$\fpt(\beta) = \wpt(\beta)$ iff $\frc(p, q) = \wrc(p, q)$;
\itemrm{(ii)}
$\ppt_k^\gh(\beta)$ is the same for all $k$ iff $\theta\big({\wrc(p,
q)}\big) = 0$.
\enditems
\endprocl

\proof Part (ii) is obvious, but part (i) needs some explanation.
One implication of (i) is also obvious from the above relations, namely,
that if $\fpt(\beta) = \wpt(\beta)$, then $\frc(p, q) = \wrc(p, q)$.
Conversely, if $\frc(p, q) = \wrc(p, q)$, then a.s.\ there cannot be more
than one infinite component. For if there were, then with positive
probability there would be neighbors $x, y$ belonging to distinct infinite
components in $\edge\setminus \{[x, y]\}$.
Call this event $A_{x, y}$.
We have $\frc(p, q)\big[ [x, y] \in \omega \mid A_{x, y}\big] =
p/[p+(1-p)q] \ne p = \wrc(p, q)\big[ [x, y] \in \omega \mid A_{x, y}\big]$,
which contradicts $\frc(p, q) = \wrc(p, q)$.

Since there cannot be more than one infinite component, the above relations
give $\fpt(\beta) = \wpt(\beta)$.
\Qed


It seems reasonable to suppose that Conjectures
\briefref g.punontriv/ and \briefref g.pcpu/ extend to random cluster
models, so that for each $q \ge 1$, there would be three phases on
nonamenable transitive graphs with one end.
In the case of the graph formed by the product of a regular tree of
sufficiently high degree and $\Z^d$, this follows from \ref
b.NewmanWu:ising/.

It is well known that $\rc(p, q)$ is stochastically increasing in $p$ for
each fixed $q$, where $\rc(p, q)$ denotes either $\frc(p, q)$ or $\wrc(p,
q)$.
Therefore, the set of $p$ for which $\theta\big({\rc(p, q)}\big) = 0$ is an
interval for each $q$. The same holds for the sets of $p$ for which the
number of infinite components is $\infty$ or 1 by the following partial
analogue of \ref t.uniq/:

\procl p.rcuniq
Let $\gh$ be a transitive unimodular graph.
Given $q \ge 1$ and $p_1 < p_2 < 1$, if there is a unique infinite
component $\rc(p_1, q)$-a.s.\ on $\gh$, then there is a unique infinite
component $\rc(p_2, q)$-a.s.
\endprocl

\proof
Given $\omega \in 2^\edge$ and $e \in \edge$, write $\omega_e$ for the
restriction of $\omega$ to $\edge\setminus \{e\}$.
A bond percolation process $(\P,\omega)$ on $\gh$ is {\bf insertion
tolerant} if $\P[e \in \omega \mid \omega_e] >0$ a.s.\ for all $e \in
\edge$.
\ref t.pu/(ii) has the following extension: If $\P$ is
any invariant ergodic percolation process on $\gh$ that is
insertion tolerant, then there is a unique infinite component $\P$-a.s.\ if
$\inf_x \P[\bp \leftrightarrow x] > 0$ (\ref B.LS:indis/).
The converse holds as well when the percolation process has positive
correlations, since then a unique infinite component implies that
$$
\P[\bp \leftrightarrow x]
\ge \P[|K(\bp)| = \infty]\P[|K(x)| = \infty] = \P[|K(\bp)|= \infty]^2
\,.
$$
We have already noted that $\rc(p, q)$ is invariant and has positive
correlations. It is easy to see that $\rc(p, q)$ is insertion tolerant, and
ergodicity is proved by \ref b.BC:covariance/ (for $\frc$) and by \ref
b.BBCK/ (for both measures).
Therefore, we may apply this extension of \ref t.pu/(ii) and its converse
to $\rc(p_1, q)$ and $\rc(p_2, q)$.
\Qed

The ergodicity needed in this proof has itself a simple proof that seems to
have been overlooked.  In fact, $\rc$ has a trivial tail $\sigma$-field on
every graph, not merely on transitive graphs.  To see this, let $B$ be any
increasing cylinder event and let $A$ be any tail event. 
Let $\gh$ be exhausted by the finite subgraphs $\gh_n$. Suppose that $n$ is
large enough that $B$ depends on the edges in $\gh_n$ only. Given $M > n$,
approximate $A$ by
a cylinder event $C$ depending only on edges in $\gh_M \setminus \gh_n$.
Let $D$ denote the event that all edges in $\gh_M \setminus \gh_n$ are closed.
Then 
$$
\rc^{\gh_n}(B) = \rc^{\gh_M}(B \mid D) \le \rc^{\gh_M}(B \mid C)
\,,
$$
because $\rc$ has positive correlations.
Letting $M \to \infty$ shows that $\rc^{\gh_n}(B) \le \frc(B \mid C)$, whence
$\rc^{\gh_n}(B) \le \frc(B \mid A)$. Now let $n \to \infty$ to conclude that
$\frc(B) \le \frc(B \mid A)$.
Since the same holds with $\neg A$ in place of $A$, this inequality is,
in fact, an equality. That is, $A$ is independent of every increasing
cylinder event, whence of every cylinder event, whence of every event. In
other words, $A$ is trivial. 
To prove tail triviality for $\wrc$, we use the same proof with $\gh_n^*$ in
place of $\gh_n$, with $D$ being the event that all edges in $\gh_M^*
\setminus \gh_n^*$ are open, and with reversed inequalities.
(A similar proof appears for different measures in \ref b.BLPS:usf/.)

There are four possible phases in Potts models that are often investigated,
i.e., four types of behavior for different values of $\beta$. We shall say
that a Potts measure at inverse temperature $\beta$ is {\bf extreme} if, as
an element of the convex set of all Potts measures on $\gh$ at inverse
temperature $\beta$, it is extreme.
The four phases are:
\smallbreak
\item{(I)}
there is a unique Potts measure (equivalently, $\ppt_k(\beta)$ does not
depend on $k$);
\item{(II)}
the free Potts measure is
extreme and there are other Potts measures;
\item{(III)}
the free Potts measure is not extreme, nor equal to the wired Potts
measure;
\item{(IV)}
the free Potts measure is equal to the wired Potts
measure and there are other Potts measures.
\smallbreak

\ref b.NewmanWu:ising/ showed the existence of three phases, namely, (I),
(II) $\cup$ (III), and (IV), each
containing an interval of parameter values of positive length, for the
$q$-state Potts model on the graph formed by the product of a
regular tree and $\Z^d$, provided that the tree has
sufficiently high degree depending on $q$.
\ref b.RS:isoD/ extended this to show that for any $q$, if $\gh$ is a
transitive graph of degree $\gdeg$ with $\isoe(\gh)/\gdeg$ sufficiently
close to 1 and with $p_u(\gh) < 1$, then there are these same
three phases in the $q$-state Potts model.
\ref b.Wu:IsingHyp/ showed similar results for a graph which is not
transitive but is similar to the hyperbolic plane.

There are some partial results for other graphs.
For natural Cayley graphs of co-compact Fuchsian groups, an uncountable
number of mutually singular Potts measures were constructed by \ref
b.SeriesSinai/.

In the following results, we use ``interval" to mean interval of positive
length.

\procl t.RC \procname{\ref B.Jonasson:RCamen/}
Let $\gh$ be a nonamenable regular graph and $q \ge 2$ be an integer.  Then
for all sufficiently large $q$, there is an interval of $p$ for which
$\frc(p, q) = \wrc(p, q)$ and there
is an interval of $p$ for which $\frc(p, q) \ne \wrc(p, q)$.
\endprocl

\noindent
[In fact, $\frc(p, q) = \wrc(p, q)$ holds for small $p$ and all $q$ since
both measures are dominated by $\P_p$. It would be interesting to show that
$\frc(p, q) = \wrc(p, q)$ can occur for large $p$ if, say, $\gh$ has one
end.]

As a consequence of this and \ref p.PtRC/(i), we obtain:

\procl t.Potts \procname{\ref B.Jonasson:RCamen/}
Let $\gh$ be a nonamenable regular graph and $q \ge 2$ be an integer.
For all sufficiently large $q$, there is an interval of $\beta$ for which
$\fpt(\beta) = \wpt(\beta)$ and there is an interval of $\beta$
for which $\fpt(\beta) \ne \wpt(\beta)$.
\endprocl

Both of the above theorems fail when $\gh$ is amenable and transitive
(\ref B.Jonasson:RCamen/).

The last result we mention also gives a characterization of amenability
among transitive graphs, but it involves an external
field. To define the Ising model with external field $h$ on a finite
graph $\gh$, modify the
energy $H(\omega)$ to be
$$
H(\omega) :=
2\sum_{e \in \edge} I_\omega(e) + 2h \sum_{x \in \vertex}
\II{\omega(x) \ne 1}
\,.
$$
Here, we take $q=2$. The corresponding probability measure on $\{-1,
1\}^\vertex$ proportional to $e^{-\beta H(\omega)}$ is denoted
$\ising^\gh(\beta, h)$. For an infinite graph $\gh$,
two limits over exhaustions $\gh_n$ are particularly
important, namely, $\ising_{\pm}^\gh(\beta, h)$, the weak${}^*$ limits of
$\ising^{\gh_n}(\beta, h)$ conditional on $\omega \restrict \bdiv \gh_n$ to
be a constant, $\pm 1$, respectively.

\procl t.ising \procname{\ref B.JonassonSteif:Ising/}
If $\gh$ is a nonamenable graph of bounded degree,
then for some $\beta$, there is an
interval of $h$ for which
$\ising_{+}^\gh(\beta, h) = \ising_{-}^\gh(\beta, h)$
and there is an
interval of $h$ for which
$\ising_{+}^\gh(\beta, h) \ne \ising_{-}^\gh(\beta, h)$.
\endprocl

As \ref b.JonassonSteif:Ising/ show, this is not true for any amenable
transitive graph.

\bsection{Percolation on Trees}{s.perctree}

As we have already mentioned,
if $T$ is a regular tree, then Bernoulli percolation produces a cluster
$K(\bp)$ that is a Galton-Watson tree (except for the first generation), so
its analysis is easy and well known.
In fact, it is not hard to find the critical value for percolation on
Galton-Watson trees:

\procl p.GWperc \procname{\ref B.Lyons:rwpt/}  Let $T$ be the family tree
of a Galton-Watson process with mean $m > 1$.  Then $p_c(T) =
1/m$ a.s.\ given nonextinction.
\endprocl

In the proof, as well as below, we write $\tu x $ for the descendant
subtree of $T$ from $x$, i.e., the tree formed from all $y \in T$ such that
the path from $\bp$ to $y$ contains $x$.

\proof
Consider Bernoulli($p$) percolation on $T$.
We claim that $K(\bp)$ has the law ({\it not\/} conditioned on $T$) of
another Galton-Watson tree having mean $mp$: Let $L$ be a random variable
whose distribution is the offspring law for $T$ and let $Y_i$ represent
i.i.d.\ Bin($1, p$) random variables that are also independent of $T$. Then
$$
\Eleft{\sum_{i=1}^L Y_i } = \Eleft{\E\bigg[ \sum_{i=1}^L Y_i \biggm|
L \bigg] } = \Eleft{\sum_{i=1}^L \E[Y_i] } =
\Eleft{\sum_{i=1}^L p } = pm \, .
$$

Hence $K(\bp)$ is finite a.s.\ if
$mp\le 1$. Since $\EBig{\Pbig{|K(\bp)| < \infty \mid T}} = \Pbig{|K(\bp)| <
\infty}$, this means that for almost every
Galton-Watson tree $T$, the component of its root is finite a.s.\ if
$mp\le 1$. In other words,
$p_c(T)\ge 1/m$ a.s.\ given nonextinction.  Similarly, the
component of the root
is infinite w.p.p.\ if $mp >1$, whence $\forall p > {1/ m}\ \
p_c(T)\le  p$ w.p.p.  It remains to show that $\P[T \hbox{ is infinite and
} p_c(T)\le
p] = 1-q$, the probability of nonextinction, for $p > {1/ m}$.
However, it is easy to see that the event $\{T \hbox{ is finite or } p_c(T)
> p\}$ is
inherited in the sense that if $T$ has this property, then so does $\tu x $
for each child $x$ of the root of $T$. It follows (e.g., see \ref
b.Lyons:book/) that
$\P[p_c(T) > p] \in\{q,1\}$.  We have already seen that it is
not equal to 1.  \Qed

Results for percolation on more
general trees depend on the following notions.
Define a {\bf cutset} of a tree $T$ to be a collection $\Pi$ of vertices
whose removal from $T$ would leave $\bp$ in a finite component.
\ref b.Lyons:rwpt/ defined the {\bf branching number} of $T$ to be
$$
\br (T)
:= \inf \left\{ \lambda > 1 \st \inf_\Pi \sum_{x \in \Pi} \lambda^{-|x|} =
0 \right\}
\,,
$$
where the infimum is over cutsets $\Pi$.
This is related to the Hausdorff dimension of the boundary of $T$: The {\bf
boundary} of $T$, denoted $\bd T$, is the set of infinite paths from $\bp$
that do not backtrack. We put a metric on $\bd T$ by letting
the distance between $\xi$ and $\eta$ be $e^{-n}$ if the number of
edges common to $\xi$ and $\eta$ is $n$.
Then $\br(T) = e^{\dim \bd T}$;
\ref b.Furstenberg/ was the first to consider $\dim \bd T$.
If $T$ is spherically symmetric (about $\bp$), meaning that $\deg x$ is a
function only of $|x|$ for $x \in T$, then $\br(T) = \gr(T)$, while in
general, we have $\br(T) \le \gr(T)$.

The following theorem was first proved (in different but equivalent
language) by \ref b.Hawkes/ for trees $T$ with bounded degree and by \ref
b.Lyons:rwpt/ in general:

\procl t.pcbr
If $T$ is any tree, then $p_c(T) = 1/\br(T)$.
\endprocl

{}From this and \ref p.GWperc/, we find that $\br(T) = m$ a.s.\ for
Galton-Watson trees with mean $m$; this was first shown (in the language of
Hausdorff dimension) by \ref b.Hawkes/.

The issue of uniqueness of infinite clusters on trees was settled in
folklore, but appeared in print for the first time by \ref b.PS:dynam/.

\procl p.num-component
For any tree $T$ and $p < 1$, the number of infinite clusters
on $T$ is $\P_p$-a.s.\ 0 or $\P_p$-a.s.\ $\infty$.
\endprocl


Similarly, one can describe the number of ends of the clusters for
percolation on trees:

\procl t.PP  \procname{\ref B.PP:critRW/}
If $T$ is any tree and $0 < p < 1$, then $\P_p$-a.s.\ either $K(\bp)$ is
finite or $K(\bp)$ has infinitely many ends.
\endprocl

In order to determine the behavior of percolation at the critical value, we
need to introduce the notion of capacity. Let $\mu$ be a probability
measure on $\bd T$. For $p < 1$, we define the $p$-{\bf energy} of $\mu$ as
$$
\en_p(\mu) :=
\int \!\! \int p^{-|\xi_1 \wedge \xi_2|} \,d\mu(\xi_1)\,d\mu(\xi_2)
\,,
$$
where $\xi_1 \wedge \xi_2$ denotes the vertex in $\xi_1 \cap \xi_2$
that is furthest from $\bp$. (If $\xi = \eta$, we interpret $p^{-|\xi_1
\wedge \xi_2|} := \infty$.)
The $p$-{\bf capacity} of $\bd T$, denoted $\CAP_{(p)}(\bd T)$,
is the reciprocal of the
minimum $\en_p(\mu)$ over all probability measures $\mu$ on $\bd T$.
If $T$ is spherically symmetric, then
$$
\CAP_{(p)}(\bd T) = \left( 1 + (1-p) \sum_{n=1}^\infty {1 \over p^{n} |\t n |}
\right)^{-1}
\,,
$$
where $\t n $ denotes the set of vertices $x$ with $|x| = n$.

The second part of
the following theorem was shown by \refbmulti{Fan:thesis,Fan:crasp} when
$T$ has bounded degree, and the full theorem by
\ref b.Lyons:rwcpt/ in general:

\procl t.perccap
If $T$ is any tree with root $\bp$ and $0 < p < 1$, then
$$
\CAP_{(p)}(\bd T) \le \theta_\bp(p) \le 2\, \CAP_{(p)}(\bd T)
\,.
$$
In particular, the probability of an infinite cluster is positive
iff $\CAP_{(p)}(\bd T) > 0$.
\endprocl

Using this result,
it is easy to construct trees $T$ for which $\theta\big(p_c(T)\big) > 0$ or for
which $\theta(p)$ is discontinuous at other $p$. Similarly, nothing like
Theorems \briefref t.merge/ or \briefref t.indist/ hold for general trees.

An extension and sharpening of \ref t.perccap/ is known for arbitrary
survival parameters. Given any survival parameters $p(e)$ on the edges $e$
of $T$, we define the corresponding energy as
$$
\en(\mu) :=
\int \!\! \int \P[\bp \leftrightarrow \xi_1 \wedge \xi_2]^{-1}
\,d\mu(\xi_1)\,d\mu(\xi_2)
$$
and define capacity as before but using this energy.
The following theorem was proved by \ref b.Lyons:rwcpt/, with the
sharpening provided by the second inequality due to \ref b.Marchal/:

\procl t.gencap
If $T$ is any tree and with any survival parameters $p(\cbuldot)$ and
corresponding capacity $\kappa := \CAP(\bd T)$, we have
$$
\kappa \le \P[\bp \leftrightarrow \infty] \le 1 - e^{-2\kappa/(1-\kappa)}
\le 2\kappa
\,.
$$
\endprocl

\ref b.HPS:dynam/ introduced a version of (bond) percolation on graphs
that evolves in time.
Given $p\in (0,1)$, the set of open edges evolves so that at any fixed time
$t\geq 0$, the distribution of this set is $\P_p$.
Let the initial distribution at time $0$ be given by $\P_p$, and let each edge
change its status (open or closed) according to a continuous-time,
stationary 2-state Markov chain, independently of all other edges.
Each edge flips (changes its value) at rate $p$ when closed and rate $1-p$
when open.
Let $\Psi_p$ denote the probability measure for this Markov process, called
{\bf dynamical percolation} with parameter $p$. This process is most
interesting for $p = p_c(\gh)$ because of the following general result:

\procl t.dynamgen \procname{\ref B.HPS:dynam/}
For any graph $\gh$, if $p > p_c(\gh)$, then $\Psi_p$-a.s.\ there is an
infinite cluster for every time $t$, while if $p < p_c(\gh)$, then
$\Psi_p$-a.s.\ there is an infinite cluster for no time $t$.
\endprocl

On trees, one can decide what happens at criticality by means of a capacity
condition (that we express for comparison via percolation instead of
capacity):

\procl t.dynam \procname{\ref B.HPS:dynam/}
Let $T$ be a tree and $0 < p < 1$. Write $\P^*$ for the probability measure
of percolation on $T$ that independently retains each edge joining $\t n-1 $
to $\t n $ with probability $p+p/n$. Then there is $\Psi_p$-a.s.\ some time
$t > 0$ at which there is an infinite cluster on $T$ iff $\P^*$-a.s.\ there
is an infinite cluster on $T$. If $T$ is spherically symmetric, then this
is equivalent to
$$
\sum_{n=1}^\infty {1 \over n p^n |\t n |} < \infty
\,.$$
\endprocl

No reasonable necessary and sufficient
condition is known so that with $\Psi_{p_c(T)}$-proba\-bi\-li\-ty 1,
there exists an infinite cluster for
{\it all\/} times $t > 0$.
However,
\ref b.PS:dynam/ have shown that when $p > p_c(T)$,
there are infinitely many infinite clusters
for all times $t$ simultaneously $\Psi_p$-a.s.
This follows from the proof of \ref p.num-component/ together with \ref
t.dynamgen/.

\bsection{The Ising Model on Trees}{s.isingtree}

We resume the notation of \ref s.ising/.
The Ising model on trees was first studied by \ref b.KKW/, who showed that
if $T$ is regular of degree $\tdeg+1$, then its critical
$\beta$ equals $\coth^{-1} \tdeg$, meaning that there is a unique Ising
measure for $\beta  > \coth^{-1} \tdeg$, but not for $\beta  < \coth^{-1}
\tdeg$ (see also \ref b.Preston/).
In other words, this is the boundary of phase (I) in the phase divisions
given in \ref s.ising/.  This calculation was extended to all trees by \ref
b.Lyons:ising/:

\procl t.tree-ising
If $T$ is any tree, then its critical $\beta$ equals $\coth^{-1} \br(T)$.
\endprocl

\ref b.BleherRZ/ showed that the critical $\beta$ for the free Ising model
on a regular tree $T$ of degree $\tdeg+1$
equals $\coth^{-1} \sqrt\tdeg$. This means that the free Ising measure is
extreme for $\beta < \coth^{-1} \sqrt\tdeg$,
but not for $\coth^{-1} \sqrt\tdeg < \beta < \coth^{-1} \tdeg$.
This is the boundary of phase (II).
A simpler proof was given by
\ref b.Ioffe:Bethe/. The result was extended to all trees by \ref b.EKPS/:

\procl t.tree-free
If $T$ is any tree, then the critical value of $\beta$ for the free Ising model
equals $\coth^{-1} \sqrt{\br(T)}$.
\endprocl

Theorem 7.7 of \ref b.Georgii:book/ shows that an Ising measure is extremal
iff it has a trivial tail, and Lemma 4.2 of \ref b.EKPS/ or Lemma 2 of \ref
b.Ioffe:ising/ shows that, for the free Ising measure,
this is equivalent to independence of $\omega(\bp)$ from the tail.
Thus, another interpretation of \ref t.tree-free/ involves asymptotic
reconstruction of $\omega(\bp)$ given $\omega(x)$ for all $x \in \t n $ as
$n \to\infty$.

We next discuss Edwards-Anderson
spin glasses. For a graph $\gh$, let $\rsign(e)$ be
independent uniform $\pm1$-valued random variables indexed by the edges $e
\in \edge$.
If $\gh$ is finite,
define the energy of a configuration $\omega \in \{1, -1\}^\vertex$ to be
$$
H(\omega) :=
2 \sum_{e \in \edge} \rsign(e) I_\omega(e)
\,.
\label e.Ham
$$
The corresponding probability measure at inverse temperature $\beta$
on $\{1, -1\}^\vertex$ is the one
proportional to $e^{-\beta H(\omega)}$, denoted $\SG^\gh(\beta)$.
Note that this measure depends on the values of $\rsign$.
We call $\P$ a {\bf spin glass measure} on an infinite graph
$\gh$ at inverse temperature $\beta$
and with interactions $\rsign(e)$
%
if $\P$ is a Markov random field and for all finite $\vertex' \subset
\vertex$ and all $\omega' \in \{1, \ldots, q\}^{\vertex'}$,
$$
\PBig{\omega \restrict \vertex' = \omega' \Bigm| \omega \restrict
\bdiv \gh' = \bd \omega'}
=
\SG^{\gh'}(\beta)\Big[\omega = \omega' \Bigm| \omega
\restrict \bdiv \gh' = \bd \omega'\Big] \,.
$$
\vbox{\noindent Define
\begineqalno
\beta^{SG}_c(\gh) :=
\sup \{ \beta \ge 0 \st &\hbox{for a.e.\ } \rsign(\cbuldot) \hbox{ and for
every } \beta' \in [0, \beta]
\cr\noalign{\vskip-3pt}&\hbox{there is a unique
spin glass measure on } \gh
\cr\noalign{\vskip-3pt}&\hbox{at inverse temperature $\beta'$
and with interactions $\rsign(\cbuldot)$} \}
\,.
\cr
\endeqalno
See \ref b.Newman:book/ for background.
}

If $T$ is a regular tree of degree $\tdeg+1$, then \ref b.CCST/ showed that
$\beta^{SG}_c(T) = \coth^{-1} \sqrt\tdeg$.
On trees, the spin glass model is
equivalent via a gauge transformation to having random
independent boundary conditions in the Ising model. Under this
transformation, let $\P_\beta$ denote the limiting Ising measure. The phase
transition defining $\beta^{SG}_c$
is equivalent to $\P_\beta$ going from not being a.s.\ extreme
to being a.s.\ extreme as $\beta$ passes $\coth^{-1} \sqrt\tdeg$.  This
calculation was extended to all trees by \ref b.PP:recursion/:

\procl t.tree-spin
If $T$ is any tree, then $\beta^{SG}_c = \coth^{-1} \sqrt{\br(T)}$.
Furthermore, there is a.s.\ more than one spin glass measure on $T$ for
every $\beta > \beta^{SG}_c$.
\endprocl

The critical cases in each of the above
three theorems can be decided based on a capacity criterion, although the
capacity for \ref t.tree-ising/ is not the usual double-integral type, but
a triple-integral type.
These capacity criteria hold for varying interaction strengths $J(e)$ as
well.
(The interaction strengths affect the Hamiltonian as in \ref e.Ham/.)
In order to state these criteria, we shall use the following notation:
For a vertex
$x \in T$ and an edge or vertex $a \in T$, write $a \le x$ if
$a$ is on the path from $\bp$ to $x$.  If $x \in \bd T$, then $a \le x$
will mean $a \in x$. Let
$$
C(x) := \prod_{e \le x} \tanh \big(J(e)\beta\big)
$$
and
$$
k(x) := \sum_{\bp \ne y \le x} C(y)^{-2}
\,.
$$

\procl t.PPising \procname{\ref B.PP:recursion/}
Let $T$ be any tree without leaves (except possibly at $\bp$).
Let $0 < \inf_{e \in \edge} J(e) \le \sup_{e \in \edge} J(e) < \infty$.
\beginitems
\itemrm{(i)} There is a unique Ising measure at inverse temperature $\beta$
iff there is a probability measure $\mu$ on $\bd T$ such that
$$
\int \!\! \int \!\! \int k(\xi_1 \wedge \xi_2 \wedge \xi_3)
\,d\mu(\xi_1)\,d\mu(\xi_2)\,d\mu(\xi_3) < \infty
\,.
$$
If $J(e) \equiv J$ and $T$ is spherically symmetric, then this is
equivalent to
$$
\sum_{n \ge 1} {1 \over [\tanh (J\beta)] ^{2n} |\t n |^2} < \infty
\,.
$$
\itemrm{(ii)} The free Ising measure at inverse temperature $\beta$
is extreme iff $\P_\beta$ is a.s.\ extreme iff there is a probability
measure $\mu$ on $\bd T$ such that
$$
\int \!\! \int k(\xi_1 \wedge \xi_2)
\,d\mu(\xi_1)\,d\mu(\xi_2) < \infty
\,.
$$
If $J(e) \equiv J$ and $T$ is spherically symmetric, then this is
equivalent to
$$
\sum_{n \ge 1} {1 \over [\tanh (J\beta)] ^{2n} |\t n |} < \infty
\,.
$$
\enditems
\endprocl

Yet another phase transition of the Ising model concerns the {\bf magnetic
susceptibility}, $\lim_{n
\to\infty} \Var\left(\sum_{x \in B(\bp, n)} \omega(x)\right) /|B(\bp, n)|$,
where $\Var$ is variance with respect to the free Ising measure.
\ref b.Matsuda/ and \ref b.Falk/ showed that the magnetic susceptibility
becomes infinite when $\beta$ passes the critical value $\coth^{-1}
\sqrt\tdeg$ if $T$ is a regular tree of degree $\tdeg + 1$.
\comment{This can be extended to general trees,
but vertices are then not all weighted equally. Compare the proof of
Theorem 1.2 of \ref b.EKPS/. But this will not appear in \ref
b.PP:recursion/.}

We turn now to models other than the Ising model.
\ref b.PemSteif:robust/ have shown that the location of a
phase transition for the $q$-state Potts model on trees with $q \ge 3$
depends on subtle aspects of the structure of the tree, and most certainly not
on $\br(T)$. However, the
location of a robust phase transition still depends only on $\br(T)$.
Here, we are using the following notion.
Given a cutset $\Pi$ of $T$, let $\Pi(\bp)$ denote the component of $\bp$
in $(T\setminus \Pi) \cup \Pi$.
For $\epsilon > 0$, let $s_\beta(\Pi, \epsilon)$ denote the distribution of
$\omega(\bp)$ with respect to the Potts measure on $\Pi(\bp)$
with inverse temperature $\beta$, interaction strengths
$$
J(e) := \cases{
\epsilon  &if $e$ has an endpoint in $\Pi$,\cr
1 & if not,\cr
}
\label e.rob
$$
and conditioned on $\omega\restrict \bdiv \Pi(\bp) \equiv 1$.
Then the critical value for a {\bf robust phase transition} is defined to
be
$$
\sup \big\{ \beta \st \all{\epsilon > 0} \inf_\Pi \|s_\beta(\Pi, \epsilon) -
1/q\|_\infty > 0 \big\}
\,,
$$
where the infimum is over all cutsets $\Pi$.
[Instead of considering arbitrarily small boundary interactions strengths
$\epsilon$, one could instead keep the interaction strengths constant and
use high temperatures at the boundaries $\Pi$.]

\procl t.Pottsrobust \procname{\ref B.PemSteif:robust/}
If $T$ is any tree with bounded degrees and $q \ge 2$, the critical value
of $\beta$ for a robust phase transition in the $q$-state Potts model on
$T$ is the unique value of $\beta$ satisfying
$$
{e^\beta +(q-1)e^{-\beta} \over e^\beta - e^{-\beta}} = \br(T)
\,.
$$
\endprocl

\noindent
In particular, the location for a robust phase transition in the Ising
model is the same as for the usual phase transition.

Lastly, we consider some continuous models on trees. Let $S^d$ denote the
$d$-dimensional unit sphere in $\R^{d+1}$. Given a finite graph $\gh$
and interaction strengths $J(e)$ for $e \in \edge$,
define the energy of $\omega \in (S^d)^\vertex$ as
$$
H(\omega) :=
\sum_{e \in \edge} H_\omega(e)
\,,
$$
where for any edge $e = [x, y]$, we write $H_\omega(e) := - J(e) \omega(x)
\cdot \omega(y)$.
The {\bf $d$-dimensional spherical measure} on $\gh$ at inverse temperature
$\beta$ is the probability measure proportional to $e^{-\beta
H(\omega)}\P(\omega)$, where $\P$ is the product measure on $(S^d)^\vertex$
with marginals on each coordinate equal to Lebesgue measure (i.e.,
normalized surface measure) on $S^d$.
When $d=1$, this is called the ``rotor" measure; when $d=2$, it is called the
``Heisenberg" measure.

For a tree $T$ and $\epsilon > 0$, let $s_\beta(\Pi, \epsilon)$ denote the
density of $\omega(\bp)$ with respect to the $d$-dimensional spherical
measure on $\Pi(\bp)$ with inverse temperature $\beta$, interaction
strengths as in \ref e.rob/, and conditioned on $\omega\restrict \bdiv
\Pi(\bp) \equiv \hat 1$, where $\hat 1$ denotes any fixed element of $S^d$.
The critical value for a {\bf robust phase transition} is defined to be
$$
\sup \big\{ \beta \st \all{\epsilon > 0} \inf_\Pi \|s_\beta(\Pi, \epsilon) -
1\|_\infty > 0 \big\}
\,.
$$

\procl t.sphere \procname{\ref B.PemSteif:robust/}
If $T$ is any tree with bounded degrees and $d \ge 1$, the critical value
of $\beta$ for a robust phase transition in the $d$-dimensional spherical
model on $T$ is the unique value of $\beta$ satisfying
$$
{\int_{-1}^1 e^{\beta r} (1-r^2)^{d/2-1}\,dr \over
\int_{-1}^1 r e^{\beta r} (1-r^2)^{d/2-1}\,dr} = \br(T)
\,.
$$
\endprocl

\bsection{The Contact Process on Trees}{s.contacttree}

The {\bf contact process} with parameter $\lambda$ on a graph $\gh$ is a
continuous-time Markov chain $\cpr_t$ on $2^\vertex$. The subset $\cpr_t
\subseteq \vertex$ is called the set of {\bf infected} (or {\bf occupied})
sites at time $t$, while $\vertex \setminus \cpr_t$ is the set of {\bf
healthy} (or {\bf vacant}) sites. Infected sites wait an exponential time
with parameter 1 and then become healthy, while a healthy site becomes
infected at a rate equal to $\lambda$ times the number of its infected
neighbors.
The measure $\P^A_\lambda$ is the measure of the above Markov chain when the
initial state is $\cpr_0 = A$. The contact process is said go {\bf extinct} if
$\P^\bp_\lambda[\all t \cpr_t \ne \emptyset] = 0$.
Otherwise, it {\bf survives}.
We make the further distinction that it {\bf survives strongly} (or {\bf
survives locally} or is {\bf recurrent}) if $\P^\bp_\lambda[\all T \texists
{t > T} \bp \in \cpr_t] > 0$, while it {\bf survives
weakly} (or {\bf globally}) if it survives but it does not survive
strongly.
It is easy to couple two copies of this Markov chain with different
parameter values so that the infected sites corresponding to the larger
value always contain the infected sites corresponding to the smaller value.
Thus, we may define
$$
\lambda_1 := \lambda_1(\gh)
:= \sup \{ \lambda \st \P^\bp_\lambda \hbox{ goes extinct} \}
=
\inf \{ \lambda \st \P^\bp_\lambda \hbox{ survives} \}
\,.
$$
We also define
$$
\lambda_2 := \lambda_2(\gh)
:=  \inf \{ \lambda \st \P^\bp_\lambda \hbox{ survives strongly} \}
\,.
$$
Thus, for any graph, we have $0 \le \lambda_1 \le \lambda_2 \le \infty$.

It is well known and
easy to show that $\lambda_1 > 1/d$ on any graph whose degrees are
bounded above by $d$: just dominate the size of the infection started from
a single site by a continuous-time branching process with mean $\lambda d$.
However, with rather small
tails in the offspring distribution, one can get $\lambda_1 = \lambda_2 =
0$ a.s.\ on Galton-Watson trees (\ref B.Pemantle:contact/).

It is significantly more difficult to study contact processes on trees than
any of the models on trees of the preceding sections. (One way to see
why this should be true is to observe that the graphical
representation of the contact process on a graph $\gh$ involves partially
oriented percolation on $\gh \times \R^+$.) Although this section is
devoted to trees, we shall briefly discuss other graphs at the end of
the section.

The first graph for which it was shown that
$0 < \lambda_1 < \lambda_2 < \infty$ was a regular tree:

\procl t.intermed
If $T$ is a regular tree of degree at least 3, then $0 < \lambda_1(T) <
\lambda_2(T) < \infty$.
\endprocl

\noindent
This was proved for trees of degree at least 4 by \ref b.Pemantle:contact/,
then for trees of degree 3 by \ref b.Liggett:contact/.  \ref
b.Stacey:contact/ gave a simpler proof of this result that extends to
certain other trees.

The following theorem describes the behavior at the critical values:

\procl t.crit
Let $T$ be a regular tree of degree $\tdeg + 1 \ge 3$ and consider the contact
process on $T$.
\beginitems
\itemrm{(i)}
There is extinction at $\lambda_1(T)$.
\itemrm{(ii)}
There is weak survival at $\lambda_2(T)$.
\enditems
\endprocl

\noindent Part (i) was shown by
\ref b.Pemantle:contact/ for $\tdeg \ge 3$ and by \ref b.MSZ:crit/ for
$\tdeg=2$.
Part (ii) was proved by \ref B.Zhang:complete/.

A basic duality property of contact processes is that for any $A, B \subset
\gh$, we have
$$
\P^A[\cpr_t \cap B \ne \emptyset] = \P^B[\cpr_t \cap A \ne \emptyset]
$$
(see \ref b.Liggett:IPS/, Theorem 1.7, Chapter VI).
When $\cpr_0 = \gh$, the distribution of $\cpr_t$ is stochastically
decreasing in time, whence it has a limit, $\bar \nu$, called the {\bf
upper stationary measure}. The {\bf lower stationary measure} is the
probability measure $\delta_\emptyset$ concentrated on the empty
configuration.
From duality, it follows that $\bar \nu = \delta_\emptyset$ iff the process
goes extinct.
One says that {\bf complete convergence} holds if for every
initial configuration $\cpr_0$, the distribution of $\cpr_t$
converges to a mixture of the lower and upper stationary measures.
In particular, when complete
convergence holds, there are no stationary measures other than the lower
and upper ones.

The argument of \ref b.Harris:set/ extends to show that if $\gh$ is
transitive, then the only automorphism-invariant extremal stationary
measures are the lower and upper ones. However, there may well be others
that are not invariant:

\procl t.weak-extremal \procname{\ref b.DurSchin:inter/, \ref
b.Zhang:complete/}
Let $T$ be a regular tree of degree at least 3. The contact process on
$T$ for $\lambda \le \lambda_1(T)$ has only one stationary measure; for
$\lambda_1(T) < \lambda \le \lambda_2(T)$, it has infinitely many extremal
stationary measures; and for $\lambda > \lambda_2(T)$, it has only two extremal
stationary measures and complete convergence holds.
\endprocl

A simpler proof of the last part of \ref t.weak-extremal/
was given by
\refbmulti{SalzSchon:second,SalzSchon:complete}.

Let $u_n(\lambda)$ be the probability that if the contact process on a
regular tree starts with one infected site at $\bp$, then a given site $x$
at distance $n$ from $\bp$ will be infected at some time. It is easy to see
that $u_{m+n}(\lambda) \ge u_m(\lambda) u_n(\lambda)$, whence
$$
\beta(\lambda) := \lim_{n \to\infty} u_n(\lambda)^{1/n}
$$
exists.
Of course, $\beta(\lambda) = 1$ when the process survives strongly. \ref
b.Liggett:Wald/ conjectured that $\beta(\lambda)
\le 1/\sqrt\tdeg$ when $\lambda \le \lambda_2(T)$. This was established by
\ref b.LalleySellke:contact/ and the equality case was determined by \ref
b.Lalley:growth-profile/:

\procl t.beta
If $T$ is a regular tree of degree $\tdeg + 1 \ge 3$, then $\beta(\lambda)
\le 1/\sqrt\tdeg$ for $\lambda \le \lambda_2(T)$, with equality
iff $\lambda = \lambda_2(T)$.
\endprocl

\ref t.beta/ implies \ref t.crit/(ii). Another
proof that $\beta(\lambda) < 1/\sqrt\tdeg$ for $\lambda < \lambda_2(T)$ was
given by \ref b.SalzSchon:complete/.
\ref t.beta/ has the following beautiful consequence for the {\bf limit
set} of $\cpr_t$, by which we mean the set of boundary points of $T$ each of
whose vertices is infected at some time. We use the same
metric on $\bd T$ as in \ref s.perctree/ for defining Hausdorff dimension
on $\bd T$.

\procl t.HD \procname{\ref b.LalleySellke:contact/, \ref
b.Lalley:growth-profile/}
If $T$ is a regular tree of degree $\tdeg + 1 \ge 3$, then the contact
process on $T$ for $\lambda_1(T) < \lambda \le \lambda_2(T)$ has a limit set on
$\bd T$ whose Hausdorff dimension is at most ${1\over2}\log \tdeg$ a.s.\ on
the event of survival, with equality iff $\lambda = \lambda_2(T)$.
\endprocl

Is it the case that $\lambda_1 = \lambda_2$ on amenable transitive graphs
and $\lambda_1 \ne \lambda_2$ on nonamenable transitive graphs? It is known
that $\lambda_1 = \lambda_2$ on the usual Cayley graphs of $\Z^d$ (\ref
B.BezGrim/).

\ref b.SalzSchon:second/ give many results for general graphs. In
particular, they prove

\procl t.extremal
Let $\gh$ be a graph of bounded degree.
\beginitems
\itemrm{(i)}
If $\lambda > \lambda_1(\Z)$, then the contact process on $\gh$ survives and
has complete convergence. In particular, $\lambda_2(\gh) \le \lambda_1(\Z)
< \infty$.
\itemrm{(ii)}
If $\gh$ is transitive and $\lambda > \lambda_2(\gh)$, then there are
exactly two extremal stationary measures.
\enditems
\endprocl

Finally, \ref b.RS:isoD/ has proved the existence of two phase transitions
on transitive graphs that are sufficiently nonamenable:

\procl t.ising2phase
If $\gh$ is a transitive graph of degree $\gdeg$ with $\isoe(\gh)/\gdeg \ge
1/\sqrt 2$, then $0 < \lambda_1(\gh) < \lambda_2(\gh) < \infty$.
\endprocl

However, \ref b.PemStac/ have exhibited nonamenable trees of bounded degree
with $\lambda_1 = \lambda_2$.


\bsection{Biased Random Walks}{s.rwl}

Given $\lambda \ge 1$, we define a nearest-neighbor
random walk on $\gh$ denoted \RWl as follows. Let $\degn x$
stand for the number of edges $[x, y]$ with $|y| = |x| - 1$.
Then the
transition probability from $x$ to an adjacent vertex $y$ is
$$
p(x, y) := \cases{\lambda /(\deg x + (\lambda - 1)\degn x) &if $|y| = |x| -
                                                1$, \cr
           1 / (\deg x + (\lambda - 1)\degn x) &otherwise. \cr}
$$
That is, from any vertex $x$,
each edge connecting $x$ to a vertex closer to $\bp$ is $ \lambda$
times more likely to be taken than any other edge incident to $x$.
(For $\lambda =1$, this
is simple random walk.) Such random walks were first
studied on trees, by \ref b.BerrettiSokal/, \ref
b.Krug/ and \ref b.LawlerSokal/.

These biased random walks are reversible and thus correspond to an electrical
network on $\gh$ (see, e.g., \ref b.DoyleSnell/, \ref b.KSK/,
Chapter IX, Section 10, or \ref b.Lyons:book/).
The conductances are given by
$$
C(x, y) := \lambda^{-(|x|\wedge |y|)}\, ,
$$
where $x$ and $y$ are adjacent vertices.

\procl t.biased \procname{\ref B.Lyons:groups/}
Let $\gh$ be a transitive graph.
If $ \lambda < \gr(\gh)$, then \RWl is transient,
while  if $ \lambda > \gr(\gh)$, then \RWl is recurrent.
Equivalently, if $ \lambda < \gr(\gh)$, then the effective conductance from
$\bp$ to infinity is positive, but not if $ \lambda > \gr(\gh)$.
\endprocl

One may also consider the {\bf rate of escape} of \RWl from $\bp$ when $
\lambda < \gr(\gh)$, i.e., $\lim_{n \to\infty} |X_n|/n$, where $X_n$ is the
location of the random walk at time $n$.  There are Cayley graphs with
$\gr(G) > 1$ but that have the surprising property that the rate of escape
of simple random walk is 0.  One example is the ``lamplighter'' group
denoted $G_1$ by \ref b.KV/.  For this example, \ref b.LPP:G1/ showed that
the rate of escape of \RWl is positive when $1 < \lambda < \gr(G_1)$. This
lack of monotonicity of behavior is quite unusual for models on transitive
graphs.  It might be that for every transitive graph $\gh$, the rate of
escape of \RWl is positive as long as $1 < \lambda < \gr(\gh)$.

The method of proof of \ref t.biased/ uses a corresponding result on trees:

\procl t.RWl-tree \procname{\ref B.Lyons:rwpt/}
Let $T$ be any tree.
If $ \lambda < \br(T)$, then \RWl is transient,
while  if $ \lambda > \br(T)$, then \RWl is recurrent.
\endprocl

The critical case in \ref t.RWl-tree/ is decided by a capacity criterion:

\procl t.RWl-cap \procname{\ref B.Lyons:rwpt/}
Let $T$ be any tree and $\lambda \ge 1$. Then \RWl is transient iff there
is a probability measure $\mu$ on $\bd T$ such that
$$
\int \!\! \int \sum_{n=0}^{|\xi_1 \wedge \xi_2|} \lambda^{n}
\,d\mu(\xi_1)\,d\mu(\xi_2) < \infty
\,.
$$
If $T$ is spherically symmetric, then this is equivalent to
$$
\sum_{n \ge 1} {\lambda^n \over |\t n |} < \infty
\,.
$$
\endprocl

\noindent Of course, when $T$ is spherically symmetric, this reduces to a
random walk on $\N$ and is well known.

An interesting model for which few results are known is that of
{\bf edge-reinforced random walk} $X_n$ ($n \ge 0$) on a graph $\gh$ with
parameter $\lambda$.
We begin with weights on all edges equal to 1.
If $X_n = x$, then $[X_n, X_{n+1}]$ is an edge incident to $x$ chosen
with probability proportional to the weights at time $n$ of the
edges incident to $x$. The weights of the edges at time $n+1$ are the same
as those at time $n$ except that the weight of $[X_n, X_{n+1}]$ is
increased by $\lambda$.
We call edge-reinforced random walk {\bf recurrent} if it returns to its
starting position infinitely often a.s.\ and {\bf transient} if it returns
to its starting position only finitely often a.s.
It seems reasonable to suppose that as $\lambda$ increases, the walk goes
from transient to recurrent as long as $\gh$ is nonamenable. The existence
and location of a phase transition was completely solved on trees by
\ref b.Pemantle:rrw/ for regular and Galton-Watson trees and by \ref
b.LP:rwre/ in general:

\procl t.errw
There is a strictly increasing
function $\lambda_E: \CO{1, \infty} \to \CO{0, \infty}$
with $\lambda_E(1) = 0$ such that if $T$ is any tree, then
edge-reinforced random walk on $T$ is transient for $\lambda <
\lambda_E\big(\br(T)\big)$ and is recurrent for $\lambda >
\lambda_E\big(\br(T)\big)$.
\endprocl

\bsection{Directions of Current Research}{s.direc}

We outline some of the themes that characterize much research in
nonamenable phase transitions and highlight some of the most important open
questions.

One contemporary theme in geometry and combinatorial group theory is the
investigation of
rough-isometry invariants (see, e.g.,
\refbmulti{Gromov:metric,Gromov:asymp}). Here, a map
$\phi : (X, d) \to (X', d')$ between metric spaces is called
a {\bf rough isometry} (or {\bf quasi-isometry}) if there are
positive constants $a$ and $b$ such that for all $x, y\in X$,
$$
a d(x, y) - a \le d'\big(\phi(x), \phi(y)\big) \le b d(x, y) + b
$$
and such that every point in $X'$ is within distance $b$ of the image of
$X$. Being roughly isometric is an equivalence relation.

In the context of graphs, we use the usual graph distance as the metric on
the vertex set. As
an example, it is easy to see that different Cayley graphs of the same
group are roughly isometric.
What properties of our models are invariant under rough isometry? For
example, in the context of Bernoulli percolation, is $p_c(\gh) < p_u(\gh)$
invariant when $\gh$ is a transitive graph? If so, \ref t.PSN/ would solve
\ref g.pcpu/ for Cayley graphs.  As another example, are critical exponents
invariant under rough isometry?  (However, they may turn out to be
the same for all nonamenable transitive graphs.) Potential-theoretic
rough-isometry invariants are known, but no nontrivial ones are known in
percolation theory.

If we specialize from rough isometries to changing generators for a
fixed group, we encounter a more refined sort of question having to do with
uniform properties: For example, it is easy to see that if $\gp$ is any finitely
generated group, then $\inf_S p_c\big(\gh(\gp, S)\big) = 0$, where the
infimum is over all finite generating sets of $\gp$ and $\gh(\gp, S)$
denotes the Cayley graph of $\gp$ with respect to $S$. But is $\sup_S
p_c\big(\gh(\gp, S)\big) < 1$? This would follow for groups of exponential
growth from \ref t.LyGp/ if it were known that $\inf_S \gr\big(\gh(\gp,
S)\big) > 1$, but this latter is an open question (see \ref b.GdlH:growth/ for
what is known about this growth problem).
No nontrivial uniform properties are known at present for, say, all
nonamenable groups.

In the other direction, rather than specializing rough isometries,
we may enlarge our equivalence classes from roughly
isometric to various classes of groups, such as nonamenable, word hyperbolic
(see \ref b.Gromov:hyper/ or \ref b.CP:hyper/), or Kazhdan [although this
last is not known to be invariant under rough isometries]. Thus, we may
search for characterizations of these classes of groups through Bernoulli
percolation or through other models, similar to \ref t.PSN/ (in combination
with the theorem of \ref b.BK:uni/ and \ref b.GKN:uni/) or Theorems
\briefref t.RC/, \briefref t.Potts/, and \briefref t.ising/.
Characterizations via invariant percolation, such as \ref t.treechar/,
would also be interesting.
As the astute reader will have observed, all of these characterizations are
of amenability only. Particularly interesting would be a characterization
of hyperbolicity. An important probabilistic characterization of Kazhdan
groups, though abstract from our point of view, is given by \ref
b.GW:Kazhdan/.

Another geometric theme concerns the appearance of spherical symmetry.
Transitive graphs are almost never spherically symmetric, i.e., it is rare
for a transitive graph to have the property that if $|x| = |y|$, then
there is an automorphism fixing $\bp$ that carries $x$ to $y$.
This lack of spherical symmetry can manifest itself in probabilistic
models. As one clear
example, $\tau_p(o, x)$ can decay to 0 as $|x| \to\infty$
in some directions while not decaying to 0 in other directions (on a given
graph); see \ref b.LS:indis/ for a
Cayley graph with this property.
What other results show the lack of spherical symmetry?
On the other hand, \ref t.biased/ has a
conclusion that holds for all spherically symmetric graphs: Here, the lack of
spherical symmetry does not affect the critical value of $\lambda$.
Are there other results where one might expect the lack of spherical
symmetry to play a role, yet where it does not? For example,
it was suggested in \ref s.rwl/
that for all transitive graphs, the rate of escape of
\RWl is positive as long as $1 < \lambda < \gr(\gh)$.

Aside from the geometrically motivated questions above, there are a plethora
of purely probabilistic questions.
The possibilities for presence or absence of
various phase transitions of random cluster and Potts
models are barely understood.
The results for contact processes that are known for trees need to be
examined for transitive graphs.  Except for branching random walks, other
interacting particle systems have barely
been investigated.

For example, we often lack
monotonicity results (such as \ref p.rcuniq/)
for processes other than Bernoulli percolation.
In fact, some such results are known to fail on
quasi-transitive graphs (see \ref b.BHW:nonmono/, for example), although there
are no known comparable failures on transitive graphs.

Finally,
some of the most basic open questions for Bernoulli percolation are:
Is $p_c(\gh) < p_u(\gh)$ when $\gh$ is a nonamenable transitive graph?
Is $p_u(\gh) < 1$ when $\gh$ is a nonamenable transitive graph with one end?
Are Theorems \briefref t.pu/(ii), \briefref t.death/, and \briefref
t.merge/ valid in the nonunimodular case?
Which transitive graphs have a unique infinite cluster at $p_u$?
What other types of phase transition are there,
such as discontinuities of $\tau_p(o, x)$
as a function of $p$ for fixed $x$?

In most situations, planar graphs are much easier to analyze due to the
availability of duality. We expect considerably faster progress for planar
graphs than for general transitive graphs.

\medbreak
\noindent {\bf Acknowledgement.}\enspace I am grateful to Itai Benjamini,
Johan Jonasson, Roberto Schonmann, Oded Schramm, Jeff Steif, the referee,
and especially
Yuval Peres for their comments. I thank Roberto and Oded for permission to
include their proofs of Theorems \briefref t.merge/ and \briefref
t.pulower/.

\ifbibnames
This order of references, which departs from the usual practice of the
Journal of Mathematical Physics, was permitted by the editor for the
convenience of the readers of the present review only.
\fi

\bibfile{\jobname}
\def\noop#1{\relax}
\input \jobname.bbl

\filbreak
\begingroup
\eightpoint\sc
\parindent=0pt\baselineskip=10pt

Department of Mathematics,
Indiana University,
Bloomington, IN 47405-5701, USA
\emailwww{rdlyons@indiana.edu}
{http://php.indiana.edu/\string~rdlyons/}

\endgroup

\bye